\documentclass[aap,MSNbibl,seceqn,dvips]{arximspdf}
\usepackage[mathscr]{eucal}

\doi{10.1214/12-AAP894} 
\volume{23}
\issue{5}
\pubyear{2013}
\firstpage{2053}
\lastpage{2098}

\makeatletter
\newcommand{\rrvert}{\vert}
\newcommand{\llvert}{\vert}

\newcommand{\IR}{\mathbb{R}}
\newcommand{\IN}{\mathbb{N}}
\newcommand{\IZ}{\mathbb{Z}}
\newcommand{\IS}{\mathbb{S}}
\newcommand{\IW}{\mathbb{W}}

\newtheorem{theorem}{Theorem}[section]
\newtheorem{corollary}{Corollary}[section]

\newproclaim{remark}{Remark}[section]

\newtheorem{lemma}{Lemma}[section]

\newproclaim{condition}{Condition}[section]

\newtheorem{proposition}{Proposition}[section]
\setattribute{abstract}{skip}{20\p@}
\makeatother

\begin{document}
\begin{frontmatter}

\title{Near critical catalyst reactant branching processes with
controlled immigration\thanksref{T1}}
\runtitle{Catalytic BP with controlled immigration}

\thankstext{T1}{Supported in part by NSF Grants
DMS-10-04418 and DMS-10-16441, the Army Research Office (W911NF-10-1-0158),
NSF Emerging Frontiers in Research and Innovation (EFRI) (Grant
CBE0736007) and the US-Israel Binational Science Foundation (Grant
2008466).}

\begin{aug}
\author[A]{\fnms{Amarjit} \snm{Budhiraja}\ead[label=e1]{budhiraj@email.unc.edu}}
\and
\author[B]{\fnms{Dominik} \snm{Reinhold}\corref{}\ead[label=e2]{dreinhold@clarku.edu}}
\runauthor{A. Budhiraja and D. Reinhold}

\affiliation{University of North Carolina at Chapel Hill and Clark University}

\address[A]{Department of Statistics\\
\quad and Operations Research\\
University of North Carolina\\
Chapel Hill, North Carolina 27599\\
USA\\
\printead{e1}}
\address[B]{Department of Mathematics\\
\quad and Computer Science\\
Clark University\\
Worcester, Massachusetts 01610\\
USA\\
\printead{e2}} 
\end{aug}

\received{\smonth{3} \syear{2012}}
\revised{\smonth{9} \syear{2012}}

%
\begin{abstract}
Near critical catalyst-reactant branching processes with controlled
immigration are studied. The reactant population evolves according to a
branching process whose branching rate is proportional to the total
mass of the catalyst. The bulk catalyst evolution is that of a
classical continuous time branching process; in addition there is a
specific form of immigration. Immigration takes place exactly when the
catalyst population falls below a certain threshold, in which case the
population is instantaneously replenished to the threshold. Such models
are motivated by problems in chemical kinetics where one wants to keep
the level of a catalyst above a certain threshold in order to maintain
a desired level of reaction activity. A diffusion limit theorem for the
scaled processes is presented, in which the catalyst limit is described
through a reflected diffusion, while the reactant limit is a diffusion
with coefficients that are functions of both the reactant and the
catalyst. Stochastic averaging principles under fast catalyst dynamics
are established. In the case where the catalyst evolves ``much faster''
than the reactant, a scaling limit, in which the reactant is described
through a one dimensional SDE with coefficients depending on the
invariant distribution of the reflected diffusion, is obtained. Proofs
rely on constrained martingale problem characterizations, Lyapunov
function constructions, moment estimates that are uniform in time and
the scaling parameter and occupation measure techniques.
\end{abstract}

%

\end{frontmatter}

\section{Introduction}\label{secintro}

This work is concerned with catalytic
branching processes that model the dynamics of catalyst-reactant
populations in which the
activity level of the reactant depends on the amount of catalyst present.
Branching processes in catalytic environments have been studied
extensively and are
motivated, for instance, by biochemical reaction networks; see \cite
{Dawson1997,Greven2009,Li2008,Kang2010} and references therein. A
typical setting consists of populations of multiple types
such that the rate of growth (depletion) of one population type is
directly affected by population sizes of other types. The simplest such model
consists of a continuous time countable state branching process
describing the evolution of the catalyst population and a second
branching process for which the branching rate is proportional to the
total mass of the catalyst population, modeling the evolution of
reactant particles.
Such processes were introduced in~\cite{Dawson1997} in the setting of
super-Brownian motions; see~\cite{Li2008}.
For classical catalyst-reactant branching processes, the
catalyst population dies out with positive probability and subsequent
to the catalyst extinction, the reactant population stays unchanged,
and therefore the
population dynamics are modeled until the time the catalyst becomes
extinct. In this work, we consider a setting where the catalyst
population is maintained above a positive threshold through a specific
form of controlled immigration. Branching process models with
immigration have also been well studied in literature; see \cite
{Athreya1972,Li2008} and references therein. However, typical
mechanisms that have been considered correspond to adding an
independent Poisson component; see, for example,~\cite{Karlin1966}.
Here, instead, we consider a model where immigration takes place only
when the population drops below a certain threshold. Roughly speaking,
we consider a sequence $\{X^{(n)}\}_{n\in\IN}$ of continuous time
branching processes, where $X^{(n)}$ starts with $n$ particles. When
the population drops below $n$, it is instantaneously restored to the
level~$n$.

There are many settings where controlled immigration models of the
above form arise naturally.
One class
of examples arises from predator-prey models in ecology,
where one may be concerned with the restoration of populations that
are close to extinction by reintroducing species when they fall
below a certain threshold.
In our work, the motivation for the study of such controlled
immigration models comes from problems in chemical reaction networks
where one wants to keep the levels of certain types of molecules above
a threshold in order to maintain a desired level of production (or
inhibition) of other chemical species in the network. Such questions
are of interest in the study of control and regulation of chemical
reaction networks. A control action where one minimally adjusts the
levels of one chemical type to keep it above a fixed threshold is one
of the
simplest regulatory mechanisms, and~the goal of this research is to
study system behavior under such mechanisms with the long-term
objective of designing optimal control policies. The specific goal of
the current work
is to derive simpler approximate and reduced models, through the theory
of diffusion approximations and stochastic averaging techniques, that
are more tractable for simulation and mathematical treatment than the
original branching process models.
In order to keep the presentation simple, we consider the setting of
one catalyst and one reactant. However, similar limit theorems can be
obtained for a more general chemical reaction network in which the
levels of some of the
chemical species are regulated in a suitable manner. Settings where
some of the chemical species act as inhibitors rather than catalysts
are also of interest and can be studied using similar techniques. These
extensions will be pursued elsewhere.

Our main goal is to establishes diffusion approximations for such
regulated catalyst-reactant systems under suitable scalings.
We consider two different scaling regimes;
in the first setting the catalyst and reactant evolve on ``comparable
timescales,'' while in the second setting the catalyst evolves ``much
faster'' than the reactant. In the former setting, the limit model is
described through a coupled system of reflected stochastic differential
equations with reflection in the space $[1,\infty)\times\IR$. The
precise result (Theorem~\ref{thmcrdiffusion}) is stated in Section
\ref{secdiffusionlimit}. Such limit theorems are of interest for
various analytic and computational reasons. It is simpler to simulate
(reflected) diffusions than branching processes, particularly for large
network settings. Analytic properties such as hitting time
probabilities and steady state behavior are more easily analyzed for
the diffusion models than for their branching process counterparts. In general,
such diffusion limits give parsimonious model representations and
provide useful qualitative insight to the underlying stochastic phenomena.

For the second scaling regime, where the catalyst evolution is much
faster, we establish
a stochastic averaging limit theorem. A key ingredient here is an
ergodicity result, which says that under a suitable ``criticality from
below'' assumption on the catalyst dynamics, the limiting catalyst
reflected diffusion admits a unique stationary distribution, which
takes an explicit form (Proposition~\ref{propstationarydistrofX}).
Characterization of
the invariant distribution is based on a variant of Echeverria's
criterion for constrained Markov processes~\cite{Kurtz1991}. Next, by
constructing suitable \textit{uniform} Lyapunov functions, we show
that the stationary distribution of the scaled catalyst branching
process converges to that of the catalyst diffusion (Theorem \ref
{thmtightnessandconvergenceofnu^n1}).
These results are then used to establish a stochastic averaging
principle that governs the dynamics of the reactant population in the
fast catalyst limit.
Proofs proceed by
developing suitable moment estimates
that are uniform in time and the scaling parameter
and by using characterization results for probability laws of reflected
diffusions through certain constrained martingale problems~\cite{Kurtz1990}.
The limit evolution of the reactant population is given through an
autonomous one-dimensional SDE with coefficients that depend on the
stationary distribution of a reflected diffusion in $[1,\infty)$. Such
model reductions are important in that they not only help in better
understanding the dynamics of the system but also help in reducing
computational costs in simulations. Indeed, since in the model
considered here the invariant distribution is explicit, the
coefficients in the one-dimensional averaged diffusion model are easily
computed, and consequently this model is significantly
easier to analyze and simulate than the original two-dimensional model.
We refer the reader to~\cite{Kang2010} and references therein for
similar results in the setting of (nonregulated) chemical reaction
networks. It will be of interest to see if similar model reductions can
be obtained for general multi-dimensional
regulated chemical-reaction networks. Key mathematical challenges will
be to identify suitable conditions for ergodicity of multi-dimensional
reflected diffusions in polyhedral domains that arise from the
regulated part of the network, and to develop uniform (in time and the
scaling parameter) moment estimates for such multi-dimensional
constrained diffusions.

We consider two different formulations of models with multiple time
scales. In Theorem~\ref{thmfastcatalystweaklimit} we consider the
setting where both catalyst and reactant processes are described
through (reflected) diffusions and the time scale parameter appears in
the coefficients of the catalyst evolution equation.
An important step here is to argue that the generator of the
two-dimensional catalyst-reactant reflected diffusion is suitably close
to the generator
of the one-dimensional averaged diffusion, for large values of the
scaling parameter. Bounds on the exponential moments of the catalyst
process, obtained in Lemma~\ref{lemunifmombdX},
play a key role in this argument.
The second formulation is considered in Theorem \ref
{fastcatalystBPaveraging}. Here, both catalyst and reactant populations evolve
according to near critical countable state branching processes, and the
branching rate in the catalyst dynamics is of higher order than that
for the reactant process. In this setting one encounters the additional
difficulty of showing that the steady state distributions of the scaled
catalyst branching process, for large values of the scaling parameter,
are suitably close to the stationary
distribution of the limiting catalyst reflected diffusion. The approach
taken here is based on characterizing the limit points of a certain
sequence of random measures
on the path space of the catalyst process and the associated reflection
process, as time and the scaling parameter together approach infinity.

The model considered in this work does not incorporate any spatial
dynamics of the two chemical species. As noted earlier in the
\hyperref[secintro]{Introduction}, in the unregulated setting, Dawson and Fleischmann \cite
{Dawson1997} considered catalyst-reactant systems, with chemical
species moving continuously in a spatial domain, given in terms of
super-Brownian motions. It will be of interest to develop analogous
continuous spatial models for the regulated catalyst-reactant systems
of the form considered in the current work. This question will be
explored in a future work.

The paper is organized as follows. We begin in Section \ref
{secdiffusionlimit} by presenting the basic limit theorem in the
setting of ``comparable time scales.'' Section \ref
{chAsymptoticBehavioroftheCatalystPopulation}
studies the time asymptotic behavior of the catalyst process under a
suitable \textit{criticality from below} assumption. Section \ref
{secDiffusionLimitoftheReactantunderFastCatalystDynamics1}
presents our main results for the multiple time scale setting. Section
\ref{secauxresults} collects some auxiliary estimates that are
needed in our proofs. Section~\ref{sec6}
proves Theorem~\ref{thmcrdiffusion}, and Section~\ref
{secproofofAsymptoticBehavioroftheCatalystPopulation} is devoted to the
proofs of Proposition~\ref{propstationarydistrofX} and Theorem
\ref{thmtightnessandconvergenceofnu^n1}. Finally, in
Section~\ref{secDiffusionLimitoftheReactantunderFastCatalystDynamics} we
present proofs of stochastic averaging principles stated
in Section~\ref{secDiffusionLimitoftheReactantunderFastCatalystDynamics1}.

\subsection{Notation}
The following notation will be used throughout this work. Denote by
$\IN$ the natural numbers, let $\IN_0:=\IN\cup\{0\}$, denote the
set of integers by $\IZ$ and let $\IR_+:=[0,\infty)$ be the set of
nonnegative real numbers. The state spaces of the scaled catalyst,
reactant, and auxiliary processes, $\hat X^{(n)},\hat
Y^{(n)}$ and $\hat Z^{(n)}$, respectively, introduced below in
(\ref{eqhatW^n}), are $\IS^{(n)}_X:=\{\frac{l}{n}|l\in
\IN_0\}\cap[1,\infty)$, $\IS^{(n)}_Y:=\{\frac{l}{n}|l\in\IN_0\}$ and
$\IS^{(n)}_Z:= \{\frac{l}{n}|l\in\IZ\}$. Let
$\IW^{(n)}:=\IS^{(n)}_X\times\IS^{(n)}_Y\times\IS^{(n)}_Z$ and
$\IW:=[1,\infty)\times\IR_+\times\IR$. Let $C^k(\IW)$ denote the space
of $k$-times continuously differentiable, real valued functions on
$\IW$, and denote by $C^k_c(\IW)$ the space of $C^k(\IW)$ functions
with compact support. Here, by a ($k$-times) differentiable function
$f$ on a set $D\subset\IR^n$ we mean a function that can be extended to
a ($k$-times) differentiable function $\tilde f$ on an open domain
$U\supset D$ such that $\tilde f$ restricted to $D$ equals $f$. Given a
metric space $S$, the space of probability measures on $S$ will be
denoted by $\mathcal{P}(S)$, the Borel $\sigma$-field on $S$ by
$\mathcal{B}(S)$, and the space of real valued, bounded, measurable
functions on $S$ by BM($S$). Let
\[
D(\IR_+\dvtx S):=\{f\dvtx \IR_+\to S|f
\mbox{ is right continuous and has left limits}\}
\]
and $D_1(\IR_+\dvtx \IR):=\{f\in D(\IR_+\dvtx \IR)|f(0)\ge1\}$, where
these $D$-spaces are endowed with the usual Skorohod topology. Let
$C(\IR_+\dvtx  \IR_+)$ be the space of continuous functions from
$\IR_+$ to $\IR_+$ endowed with the local uniform topology. We say a
sequence $\{\xi_n\}_{n\in\IN}$ of random variables with values in some
Polish space $\mathcal E$ is tight if the corresponding probability
laws are a tight sequence in $\mathcal{P}(\mathcal{E})$. For a function
$\xi\dvtx \IR_+\to\IR^n$, let the jump at time $t$ be defined as
$\Delta\xi_t:=\xi_t-\xi_{t-}$, $t>0$, and $\Delta\xi_0:=0$. For a
function $f\dvtx \IR_+\to\IR$ and $t\ge0$, let $|f|_{*,t}:=\sup_{s\le
t}|f(s)|$. For two semimartingales $\xi$ and $\zeta$, the quadratic
covariation (or bracket process) and predictable (or conditional)
quadratic covariation are denoted by $\{[\xi,\zeta]_t\}_{t\in\IR_+}$
and $\{ \langle\xi,\zeta\rangle_t\}_{t\in\IR_+}$, respectively; their
definition will be recalled in Section~\ref{secauxresults}.

\section{Diffusion limit under comparable timescales}\label{secdiffusionlimit}

Consider\vspace*{1pt} a sequence of pairs of continuous time, countable state Markov
branching
processes $(X^{(n)},Y^{(n)})$, where $X^{(n)}$ and $Y^{(n)}$ represent
the number of catalyst and reactant particles, respectively.
The dynamics are described as follows.
Each of the $X^{(n)}_t$ particles alive at time $t$ has an
exponentially distributed lifetime with parameter $\lambda^{(n)}_1$ (mean
lifetime $1/\lambda^{(n)}_1$). When it
dies, each such particle gives rise to a\vadjust{\goodbreak} number of offspring, according
to the
offspring distribution $\mu_1^{(n)}(\cdot)$. Additionally, if the
catalyst population drops below $n$, it is instantaneously replenished
back to the level $n$ (\textit{controlled immigration}). The branching
rate of the reactant process $Y^{(n)}$ is of the order of the current
total mass of the catalyst population, that is, $X^{(n)}/n$, and we
denote the offspring distribution of $Y^{(n)}$ by $\mu_2^{(n)}(\cdot
)$. A precise definition of the pair $(X^{(n)},Y^{(n)})$ will be given below.
We are interested in the study of asymptotic behavior of $(X^{(n)},
Y^{(n)})$, under suitable scaling, as $n\to\infty$.

To facilitate some weak convergence arguments, we will consider an auxiliary
sequence of processes $Z^{(n)}$ that ``shadow'' $X^{(n)}$ in the
following manner. The process $Z^{(n)}$ will be a $\IZ$ valued pure
jump process whose jump instances and sizes are the same as that of
$X^{(n)}$ away from the boundary $\{n\}$, whereas when $X^{(n)}$ is at
the boundary, $Z^{(n)}$ has a negative jump of size 1 whenever there is
immigration of a catalyst particle into the system. This description is
made precise through the infinitesimal generator given in (\ref
{eqgeneratorofwholeprocess}). The process $\hat Z^{(n)}$
will not appear in the statements of the results; nevertheless it plays
an important role in our proofs.

We now give a precise description of the various processes and the
scaling that is considered.
Roughly speaking, time is accelerated by a factor of $n$, and mass is
scaled down by a factor of $n$. Define RCLL processes
%
\begin{equation}
\label{eqhatW^n} \hat{\mathbf W}^{(n)}_t:=
\bigl(\hat X_t^{(n)},\hat Y_t^{(n)},\hat
Z_t^{(n)} \bigr):= \biggl(\frac{X_{nt}^{(n)}}{n},
\frac
{Y_{nt}^{(n)}}{n},\frac{Z_{nt}^{(n)}}{n} \biggr),\qquad t\in\IR_+,
\end{equation}
and let $\hat{\mathbf W}^{(n)}_0 = (x^{(n)}_0,y^{(n)}_0,z^{(n)}_0)\in\IW^{(n)}$, where $(nx^{(n)}_0,ny^{(n)}_0)$ is the initial number of
catalyst and reactant particles and $z^{(n)}_0=x^{(n)}_0$.
Then
$\{\hat{\mathbf W}_t^{(n)}\}_{t\in\IR_+}$ is characterized as the $\IW^{(n)}$ valued
Markov process with sample paths in
$D(\IR_+\dvtx \IW^{(n)})$, starting at $\hat{\mathbf
W}^{(n)}_0=(x^{(n)}_0,y^{(n)}_0,z^{(n)}_0)$, and having infinitesimal
generator $\mathcal{\hat A}^{(n)}$ given as
%
\begin{eqnarray}
\label{eqgeneratorofwholeprocess}\quad
\mathcal{\hat A}^{(n)}\phi(\mathbf w)
&=&\lambda_1^{(n)} n^2 x\sum
_{k=0}^\infty \biggl[\phi \biggl(1\vee x+
\frac{k-1}{n},y,z+\frac
{k-1}{n} \biggr)-\phi(\mathbf w) \biggr]
\mu_1^{(n)}(k)
\nonumber\\[-8pt]\\[-8pt]
&&{} +\lambda_2^{(n)} n^2 x y \sum
_{k=0}^\infty \biggl[\phi \biggl(x,y+
\frac{k-1}{n},z \biggr)-\phi(\mathbf w) \biggr]\mu_2^{(n)}(k),\nonumber
\end{eqnarray}
where $\mathbf w=(x,y,z)\in\IW^{(n)}$ and $\phi\in\operatorname
{BM}(\IW)$.
From the definition of the generator we see that, for each $k \ge0$,
given $\hat{\mathbf W}^{(n)}_t=(x,y,z)\in\IW^{(n)}$, the process jumps
to $(x,y+\frac{k-1}{n},z)$ with rate $\lambda^{(n)}_2 n^2xy\mu^{(n)}_2(k)$ and to $(x+\frac{k-1}{n},y,z+\frac{k-1}{n})$ with rate
$\lambda^{(n)}_1 n^2x\mu^{(n)}_1(k)$,
except when $k=0$ and $x=1$, in which case the latter jump is to
$(x,y,z+\frac{k-1}{n})$ with rate $\lambda^{(n)}_1 n^2\mu^{(n)}_1(0)$.
This property of the generator at $x=1$ accounts for the instantaneous
replenishment of the (unscaled) catalyst population to level $n$,
whenever the catalyst drops below~$n$.

For $i=1,2$, let
\[
m_i^{(n)}:=\sum_{k=0}^\infty
k\mu_i^{(n)}(k) \quad\mbox{and}\quad \alpha_i^{(n)}=
\sum_{k=0}^\infty(k-1)^2
\mu_i^{(n)}(k).
\]
We make the following basic assumption on the parameters of the
branching rates and offspring distributions as well as on the initial
configurations of the catalyst and reactant populations:
%
\begin{condition}\label{conditioncr}
\textup{(i)} For $i=1,2$ and for $n\in\IN$, $\alpha_i^{(n)},\lambda_i^{(n)}\in(0,\infty)$
and
$m_i^{(n)}=1+\frac{c_i^{(n)}}{n}, c_i^{(n)}\in(-n,\infty)$.

\mbox{}\hphantom{\textup{i}}\textup{(ii)} For $i=1,2$, as $n\to\infty$, $c_i^{(n)}\to c_i\in\IR
$, $\alpha^{(n)}_i\to\alpha_i\in(0,\infty)$ and $\lambda^{(n)}_i\to\lambda_i\in(0,\infty)$.

\textup{(iii)} For $i=1,2$ and for every $\varepsilon\in(0,\infty)$,
\[
\lim_{n\to\infty}\sum_{l:l>\varepsilon\sqrt{n}}\bigl(l-m^{(n)}_i
\bigr)^2\mu_i^{(n)}(l)= 0.
\]

\mbox{}\hspace*{1pt}\textup{(iv)} As $n\to\infty$, $ (x_0^{(n)},y_0^{(n)})\to(x_0,y_0)\in
[1,\infty)\times\IR_+$.
\end{condition}
Condition~\ref{conditioncr} and the form of the generator in (\ref
{eqgeneratorofwholeprocess}) ensure that the scaled catalyst and
reactant processes transition on comparable time scales, namely
$\mathcal O(n^2)$. In order to state the limit theorem for $(\hat
X^{(n)}, \hat Y^{(n)})$, we need some notation and definitions
associated with the one-dimensional Skorohod map with reflection at~1.
Let $\Gamma\dvtx D_1(\IR_+\dvtx \IR)\rightarrow D(\IR_+\dvtx
[1,\infty))$ be defined as
%
\begin{equation}
\label{eqSkorohodmap} \Gamma(\psi) (t):=\bigl(\psi(t)+1\bigr)-
\inf_{0\leq s\leq t}\bigl\{\psi(s)\wedge 1\bigr\}\qquad \mbox{for } \psi\in D(\IR_+\dvtx
\IR).
\end{equation}
The function $\Gamma$, known as Skorohod map, can be characterized as
follows; see, for example, Appendix B in~\cite{BeWi1} and references therein:
if $\psi,\phi,\eta^*\in D(\IR_+\dvtx \IR)$ are such that (i) $\psi
(0)\ge1$, (ii) $\phi=\psi+\eta^*$, (iii) $\phi\geq1$, (iv) $\eta^*$ is nondecreasing, $\int_{[0,\infty)} 1_{\{\phi(s)\neq1\}}\,d\eta^*(s)=0$, and $\eta^*(0)=0$, then $\phi=\Gamma(\psi)$ and $\eta^*=\phi-\psi$. The process $\eta^*$ can be regarded as the
reflection term that is applied to the original trajectory $\psi$ to
produce a trajectory $\phi$ that is constrained to $[1,\infty)$.
From the definition of the Skorohod map and using the triangle
inequality, we get the following Lipschitz property: for $\psi,\tilde
\psi\in D_1(\IR_+\dvtx \IR)$,
%
\begin{equation}
\label{eqSkorohodLipschitz} \sup_{s\leq t}\bigl|\Gamma(\psi) (s)-\Gamma(\tilde
\psi) (s)\bigr|\leq2\sup_{s\leq t}\bigl|\psi(s)-\tilde\psi(s)\bigr|.
\end{equation}

The diffusion limit of $(\hat X^{(n)},\hat Y^{(n)})$ will be the
process $(X,Y)$, starting at $(x_0,y_0)$, which is given through a
system of stochastic integral equations as in the following proposition.
%
\begin{proposition}\label{propuniqunessofX,Y}
Let $(\bar\Omega, \bar{\mathcal{F}},\bar P, \{\bar{\mathcal
{F}}_t\})$ be a filtered probability space on which are given
independent standard $\{\bar{\mathcal{F}}_t\}$ Brownian\vadjust{\goodbreak} motions $B^X$
and $B^Y$.
Let $X_0, Y_0$ be square integrable $\bar{\mathcal{F}}_0$ measurable
random variables with values in $[1,\infty)$ and $\IR_+$, respectively.
Then the following system of stochastic integral equations has a unique
strong solution:
%
\begin{eqnarray}
\label{eqX}
X_t&=&\Gamma \biggl(X_0+ \int_0^\cdot c_1\lambda_1X_s \,ds +
\int _0^\cdot\sqrt{\alpha_1 \lambda_1 X_s}\,dB^X_s \biggr) (t),
\\
\label{eqY}
Y_t&=&Y_0+\int_0^t c_2\lambda_2 X_s Y_s \,ds +\int
_0^t\sqrt{\alpha_2 \lambda_2 X_s Y_s}\,dB^Y_s,
\\
\label{eqeta} \eta_t&=&X_t-X_0-\int_0^t c_1\lambda_1X_s \,ds -\int _0^t
\sqrt{\alpha_1 \lambda_1X_s}\,dB^X_s,
\end{eqnarray}
where $\Gamma$ is the Skorohod map defined in (\ref{eqSkorohodmap}).
\end{proposition}

In the above proposition, by a strong solution of (\ref{eqX})--(\ref
{eqeta}), we mean an $\bar{\mathcal{F}}$-adapted continuous process
$(X,Y,\eta)$ with values in $[1,\infty)\times\IR_+\times\IR_+$
that satisfies (\ref{eqX})--(\ref{eqeta}).
The following is the main result of this section.

\begin{theorem}\label{thmcrdiffusion}
Suppose Condition~\ref{conditioncr} holds. The process $(\hat
X^{(n)},\hat Y^{(n)})$ converges weakly in $D(\IR_+\dvtx
[1,\infty)\times\IR_+)$ to the process $(X,Y)$ given in
Proposition~\ref{propuniqunessofX,Y} with $(X_0,Y_0)=(x_0,y_0)$.
\end{theorem}

Proposition
\ref{propuniqunessofX,Y} follows by standard arguments, so its
proof is relegated to the \hyperref[app]{Appendix}. Theorem~\ref{thmcrdiffusion} will
be proved in Section~\ref{sec6}.

\section{Asymptotic behavior of the catalyst population}
\label{chAsymptoticBehavioroftheCatalystPopulation}
Stochastic averaging results in this work rely on understanding the
time asymptotic behavior of the catalyst process. Such behavior, of
course, is also of independent interest.
We begin with the following result on the stationary distribution of
$X$, where $X$ is the reflected diffusion from Proposition \ref
{propuniqunessofX,Y}, approximating the catalyst dynamics \mbox{(Theorem
\ref{thmcrdiffusion})}.
The proof uses an extension of the Echeverria criterion for stationary
distributions of diffusions to the setting of constrained diffusions;
see Section~\ref{secproofstationarydistrofX}.
We will make the following additional assumption. Recall the constants
$c^{(n)}_1\in(-n,\infty)$ and $c_1\in\IR$ introduced in Condition
\ref{conditioncr}.
%
\begin{condition}\label{conditionc1}
For all $n\in\IN$, $c^{(n)}_1<0$ and $c_1<0$.
\end{condition}

\begin{proposition}\label{propstationarydistrofX}
Suppose Condition~\ref{conditionc1} holds. The process $X$ defined
through (\ref{eqX}) has a unique stationary distribution, $\nu_1$,
which has density
%
\begin{equation}
\label{eqdensityofnu1} p(x):= %
\cases{\displaystyle \frac{\theta}{x}\exp\biggl(2
\frac{c_1}{\alpha_1}x\biggr), &\quad if $x\geq1$,
\vspace*{2pt}\cr
0, &\quad if $x<1$, } %
\end{equation}
where $\theta:= (\int_1^\infty(\frac{1}{x}\exp(2\frac
{c_1}{\alpha_1}x))\,dx )^{-1}$.\vadjust{\goodbreak}
\end{proposition}

The following result shows that the time asymptotic behavior of the
catalyst population is well approximated by that of its diffusion
approximation given through (\ref{eqX}).
We make the following additional assumption on the moment generating
function of the offspring distribution, which will allow us to
construct certain ``uniform Lyapunov functions''
that play a key role in the analysis;
see Theorem~\ref{thmintegratingtillreturntime} and the
function $\hat V^{(n)}$ defined in (\ref{lyapfn}).

\begin{condition}\label{conditionmoment}
For some $\bar\delta>0$,
%
\begin{equation}
\label{eqoffsprdistrhasboundedmgf} \sup_{n\in\IN} \sum
_{k=0}^\infty e^{\bar\delta k}\mu^{(n)}_1(k)<
\infty.
\end{equation}
\end{condition}

\begin{theorem}\label{thmtightnessandconvergenceofnu^n1}
Suppose Conditions~\ref{conditioncr},~\ref{conditionc1} and \ref
{conditionmoment} hold. Then, for each $n\in\IN$, the process $\hat
X^{(n)}$ has a unique stationary distribution $\nu^{(n)}_1$, and the
family $\{\nu^{(n)}_1\}_{n\in\IN}$ is tight.
As $n\to\infty$, $\nu^{(n)}_1$ converges weakly to $\nu_1$.
\end{theorem}

Proposition~\ref{propstationarydistrofX} and Theorem \ref
{thmtightnessandconvergenceofnu^n1} will be proved in Section
\ref{secproofofAsymptoticBehavioroftheCatalystPopulation}.

\section{Diffusion limit of the reactant under fast catalyst
dynamics}\label{secDiffusionLimitoftheReactantunderFastCatalystDynamics1}

As noted in Section~\ref{secdiffusionlimit},
the catalyst and reactant populations whose scaled evolution is
described through (\ref{eqgeneratorofwholeprocess}) transition on
comparable time scales. In situations in which the catalyst evolves
``much faster'' than the reactant, one can hope to find a simplified
model that captures the dynamics of the reactant population in a more
economical fashion. One would expect that the reactant population can
be approximated by a diffusion whose coefficients depend on the
catalyst only through the catalyst's stationary distribution. Indeed,
we will show that the (scaled) reactant population can be approximated by
the solution of
%
\begin{equation}
\label{eqcheckY} \check Y_t=\check Y_0+\int
_0^t c_2\lambda_2
m_X \check Y_s \,ds+\int_0^t
\sqrt{\alpha_2\lambda_2 m_X\check
Y_s}\,dB_s,\qquad \check Y_0=y_0,
\end{equation}
where
$m_X=\int_1^\infty x \nu_1(dx) = - \frac{\alpha_1 \theta}{2c_1}
\exp(2c_1/\alpha_1)$.\vspace*{2pt}

Such model reductions (see~\cite{Kang2010} and references therein for
the setting of chemical reaction networks) not only help in better
understanding the dynamics of the system but also help in reducing
computational costs in simulations.
In this section we will consider such stochastic averaging results in
two model settings.
First, in Section~\ref{secStochAveraginginDiffusionSetting}, we
consider the simpler setting where the population mass evolutions are
described through (reflected) stochastic integral equations and a
scaling parameter in the coefficients of the model distinguishes the
time scales of the two processes. In Section \ref
{secStochasticAveragingforScaledBranchingProcesses} we will consider a setting
which captures the underlying physical dynamics more accurately in the
sense that the mass processes are described in terms of continuous time
branching processes, rather than diffusions.\looseness=1

\subsection{Stochastic averaging in a diffusion setting}
\label{secStochAveraginginDiffusionSetting}
In this section we consider the setting where the catalyst and reactant
populations evolve according to (reflected) diffusions similar to $X$
and $Y$ from Proposition~\ref{propuniqunessofX,Y}, but where the
evolution of the catalyst is accelerated by a factor of $a_n$ such that
$a_n\uparrow\infty$ as $n\uparrow\infty$ (i.e., drift and
diffusion coefficients are scaled by $a_n$).
More precisely, we consider a system of catalyst and reactant
populations that are given as solutions of
the following system of stochastic integral equations: for $t\ge0$,
\begin{eqnarray*}
\check X^{(n)}_t&=&\Gamma \biggl(\check X^{(n)}_0
+ \int_0^\cdot a_n
c_1\lambda_1\check X^{(n)}_s \,ds
+ \int_0^\cdot\sqrt{a_n
\alpha_1\lambda_1 \check X^{(n)}_s}\,dB^X_s
\biggr) (t),
\\
\check Y^{(n)}_t&=&\check Y^{(n)}_0+
\int_0^t c_2\lambda_2
\check X^{(n)}_s \check Y^{(n)}_s \,ds +
\int_0^t\sqrt{\alpha_2
\lambda_2 \check X^{(n)}_s \check
Y^{(n)}_s}\,dB^Y_s,
\end{eqnarray*}
where $(\check X^{(n)}_0,\check Y^{(n)}_0)=(x_0,y_0)$, $c_1, c_2 \in
\IR$, $\alpha_i,\lambda_i \in(0,\infty)$, $B^X$ and $B^Y$ are
independent standard Brownian motions, and $\Gamma$ is the Skorohod
map described above Proposition~\ref{propuniqunessofX,Y}.

The following result says that if $c_1<0$, then the reactant population
process~$\check Y^{(n)}$, which is given through a coupled
two-dimensional system,
can be well approximated by the one-dimensional diffusion $\check Y$ in
(\ref{eqcheckY}), whose coefficients are given in terms of the
stationary distribution of the catalyst process.

\begin{theorem}\label{thmfastcatalystweaklimit}
Suppose Condition~\ref{conditionc1} holds. The process $\check
Y^{(n)} $ converges weakly in $C(\IR_+\dvtx  \IR_+)$ to the process
$\check Y$.
\end{theorem}
The proof of Theorem~\ref{thmfastcatalystweaklimit} is given in
Section~\ref{secDiffusionLimitoftheReactantunderFastCatalystDynamics}.

\subsection{Stochastic averaging for scaled branching processes}
\label{secStochasticAveragingforScaledBranchingProcesses}

We now consider stochastic averaging for
the setting where the catalyst and reactant populations are described
through branching processes.
Consider catalyst and reactant populations evolving according to the
branching processes introduced in Section~\ref{secdiffusionlimit},
but where the catalyst evolution is sped up by a factor of $a_n$ such
that $a_n\uparrow\infty$ monotonically as $n\uparrow\infty$. That
is, we consider a sequence of catalyst populations $\tilde
X^{(n)}_t:=X^{(n)}_{a_nt}$, $t\ge0$, where $X^{(n)}$ are the branching
processes introduced in Section~\ref{secdiffusionlimit}. The
reactant population evolves according to a branching process, $\tilde
Y^{(n)}$, whose branching rate, as before, is of the order of the
current total mass of the catalyst population,
$\tilde X^{(n)}/n$. The infinitesimal generator
$\check{\mathcal{G}}^{(n)}$ of the scaled process
\[
\bigl(\check X^{(n)}_t,\check Y^{(n)}_t
\bigr):= \biggl(\frac
{1}{n}\tilde X^{(n)}_{nt},
\frac{1}{n}\tilde Y^{(n)}_{nt} \biggr),\qquad t\geq0,
\]
is given as
%
\begin{eqnarray}
\label{eqcheckG^n} \quad\check{\mathcal{G}}^{(n)}\phi(x,y)
&=&\lambda_1^{(n)} n^2 a_n
x\sum_{k=0}^\infty \biggl[\phi \biggl( 1\vee
\biggl(x+\frac{k-1}{n}\biggr),y \biggr)-\phi(x,y) \biggr]\mu_1^{(n)}(k)
\nonumber\\[-8pt]\\[-8pt]
&&{} +\lambda_2^{(n)} n^2 x y \sum
_{k=0}^\infty \biggl[\phi \biggl(x,y+\frac{k-1}{n}
\biggr)-\phi(x,y) \biggr]\mu_2^{(n)}(k),\nonumber
\end{eqnarray}
where $(x,y)\in\IS^{(n)}_X\times\IS^{(n)}_Y$ and $\phi\in
\operatorname{BM}([1,\infty)\times\IR_+)$.

We note that a key difference between the generators $\check{\mathcal{G}}^{(n)}$ above and $ \mathcal{\hat A}^{(n)}$ in (\ref
{eqgeneratorofwholeprocess}) is the extra factor of $a_n$ in the first
term of
(\ref{eqcheckG^n}), which says that, for large $n$, the catalyst
dynamics are much faster than that of the reactant.

We will show in Theorem~\ref{fastcatalystBPaveraging} that the
reactant population process $\check Y^{(n)}$ can be well approximated
by the one-dimensional diffusion $\check Y$ in (\ref{eqcheckY}).
Once again, the result provides a model reduction that is potentially
useful for simulations and also for a general qualitative understanding
of reactant dynamics near criticality.
%
\begin{theorem}\label{fastcatalystBPaveraging}
Suppose Conditions~\ref{conditioncr},~\ref{conditionc1} and \ref
{conditionmoment} hold. Then,
as $n\to\infty$, $\check Y^{(n)} $ converges weakly in $D(\IR_+\dvtx  \IR_+)$ to the process $\check Y$.
\end{theorem}
We will prove the above theorem in Section \ref
{secDiffusionLimitoftheReactantunderFastCatalystDynamics}.

\section{Auxiliary results} \label{secauxresults}

In this section we collect several auxiliary results, which will be
used in the proofs of our main results.
Recall that
the quadratic covariation (or bracket process) of two semimartingales
$\xi$ and $\zeta$ is the process $\{[\xi,\zeta]_t\}_{t\in\IR_+}$
defined by
\[
[\xi,\zeta]_t:=\xi_t\zeta_t-\int
_0^t \xi_{s-} \,d
\zeta_s-\int_0^t
\zeta_{s-} \,d\xi_s,\qquad t\geq0,
\]
where $\xi_{0-}:=0, \zeta_{0-}:=0$. The predictable quadratic
covariation of $\xi$ and $\zeta$ is the unique predictable process $\{
\langle\xi, \zeta\rangle_t\}_{t\in\IR_+}$ such that
$\{[\xi,\zeta]_t-\langle\xi,\zeta\rangle_t\}_{t\in\IR_+}$
is a local martingale. If $\xi=\zeta$, then $[\xi]\equiv[\xi,\xi
]$ and $\langle\xi\rangle\equiv\langle\xi,\xi\rangle$ are,
respectively, the quadratic and predictable quadratic variation
processes of $\xi$.

For $\mathbf x=(x_1,x_2,x_3)\in\IW$, let $\phi_i(\mathbf x)=x_i$, $i=1,2,3$,
and $h:=\phi_1-\phi_3$. Note that for a locally bounded measurable
function $f$ on $\IW$
%
\begin{equation}
\label{locmart} M_t^{(n)}(f):=f\bigl(\hat{\mathbf
W}_t^{(n)}\bigr)-f\bigl(\hat{\mathbf W}_0^{(n)}
\bigr)-\int_0^t\mathcal{\hat
A}^{(n)}f\bigl(\hat{\mathbf W}_s^{(n)}\bigr)\,ds,\qquad t
\ge0,
\end{equation}
is a local martingale with respect to the filtration $\sigma(\hat
{\mathbf W}^{(n)}_s\dvtx  s\le t)$. For the rest of the paper, we suppress the
filtration, and simply refer to $M^{(n)}(f)$ as a local martingale.

Let
%
\begin{equation}
\label{hateta^{n}} \hat\eta_t^{(n)}:=
\lambda_1^{(n)}n\mu_1^{(n)}(0)\int
_0^t 1_{\{\hat
X_s^{(n)}=1\}}\,ds.
\end{equation}
This process will play the role of the reflection term in the dynamics
of the catalyst, arising from the controlled immigration.
The following tightness result will be used in the weak convergence proofs.

\begin{proposition}\label{lemtighnessofcatalyst-reactant}
Suppose Conditions~\ref{conditioncr} and~\ref{conditionc1} hold.
Then the family $\{(\hat X^{(n)},\hat Y^{(n)},\hat\eta^{(n)})\}_{n\in
\IN}$ is tight in $D(\IR_+\dvtx [1,\infty)\times\IR_+\times\IR_+)$.
If additionally Condition~\ref{conditionmoment} holds, then the
family $\{(\hat X^{(n)}_{s+\cdot},\hat\eta^{(n)}_{s+\cdot}-\hat
\eta^{(n)}_s)\}_{n\in\IN, s\in\IR_+}$ is tight in $D(\IR_+\dvtx [1,\infty)\times\IR_+)$.
\end{proposition}

The proof of Proposition~\ref{lemtighnessofcatalyst-reactant} will
be based on the following results. Lem\-ma~\ref
{lemrepresentationofcrprocesses} below gives some useful
representations for the catalyst and
reactant processes.
Lemmas~\ref{lemEVofsquaredsuprema}--\ref
{lemEVofsquaredsupremaofshiftedprocesses} and Corollary \ref
{corEVofsquaredsupremaofshiftedhatX^n} provide moment bounds that are
useful for
arguing tightness. Proofs of these results
are given in Section~\ref{secproofsofauxresults}.

\begin{lemma}\label{lemrepresentationofcrprocesses}
Suppose Condition~\ref{conditioncr}\textup{(i)} holds. The process $(\hat
X^{(n)}, \hat Y^{(n)})$ can be represented as
%
\begin{eqnarray}
\label{eq5.3} \hat X_t^{(n)}&=&\hat X^{(n)}_0+c_1^{(n)}
\lambda_1^{(n)}\int_0^t
\hat X_s^{(n)}\,ds+M_t^{(n)}(
\phi_1)+\hat\eta_t^{(n)}
\nonumber\\[-8pt]\\[-8pt]
&=&\Gamma \biggl(\hat X_0^{(n)}+c_1^{(n)}
\lambda_1^{(n)}\int_0^\cdot
\hat X_s^{(n)}\,ds+M_\cdot^{(n)}(
\phi_1) \biggr) (t)\nonumber
\end{eqnarray}
and
%
\begin{equation}
\label{eqhatY^nb^n2} \hat Y^{(n)}_t=
\hat Y^{(n)}_0+c_2^{(n)}
\lambda_2^{(n)} \int_0^t
\hat X_t^{(n)}\hat Y_t^{(n)}\,ds+M^{(n)}_t(
\phi_2).
\end{equation}
Moreover, for $t \ge0$,
%
\begin{eqnarray}
\label{eqpqvM^nphi1} \bigl\langle M^{(n)}(
\phi_1)\bigr\rangle_t&=&\lambda_1^{(n)}
\alpha_1^{(n)}\int_0^t
\hat X_s^{(n)}\,ds-\lambda^{(n)}_1
\mu^{(n)}_1(0)\int_0^t
1_{\{\hat
X^{(n)}_s=1\}}\,ds
\nonumber\\[-8pt]\\[-8pt]
&\le&\lambda_1^{(n)}\alpha_1^{(n)}
\int_0^t\hat X_s^{(n)}\,ds\nonumber
\end{eqnarray}
and
%
\begin{equation}
\label{eq5.6} \bigl\langle M^{(n)}(\phi_2)\bigr
\rangle_t=\lambda^{(n)}_2\alpha_2^{(n)}
\int_0^t\hat X^{(n)}_s
\hat Y^{(n)}_s \,ds.
\end{equation}
\end{lemma}

Let
%
\begin{equation}
\label{eqhatN^ndef} \hat N_t^{(n)}:=\hat
X_0^{(n)}+c_1^{(n)}
\lambda_1^{(n)}\int_0^t
\hat X_s^{(n)}\,ds+M^{(n)}_t(
\phi_1).
\end{equation}
Then we have the following second moment estimate.
%
\begin{lemma}\label{lemEVofsquaredsuprema}
Suppose Conditions~\ref{conditioncr}\textup{(i)} and \textup{(ii)} hold. Then there is
a $K\in(0,\infty)$ such that for all $n\in\IN$ and $T\geq0$,
%
\begin{equation}
\label{eqEVofsquaredsuprema}
 E \Bigl(\sup_{t\leq T} \bigl(\bigl(\hat
X_t^{(n)}\bigr)^2+ \bigl(M^{(n)}_t(
\phi_1)\bigr)^2+ \bigl(\hat N^{(n)}_t
\bigr)^2+\bigl(\hat\eta_t^{(n)}
\bigr)^2 \bigr) \Bigr)\leq\exp\bigl(K T^2\bigr)
\bigl(x^{(n)}_0\bigr)^2\hspace*{-30pt}
\end{equation}
and for each $k\in\IN$,
%
\begin{equation}
\label{eqestimateonhatY^n} E \Bigl(\sup_{t\leq T} \bigl(\hat
Y_{\sigma_k^{(n)}\wedge
t}^{(n)} \bigr)^2 \Bigr)\leq\exp\bigl(K
T^2 k^2\bigr) \bigl(y^{(n)}_0
\bigr)^2,
\end{equation}
where $\sigma^{(n)}_k:=\inf\{t>0\dvtx \hat X^{(n)}_t\ge k\}$.
\end{lemma}

In order to study properties of invariant measures of $\hat X^{(n)}$,
it will be convenient to allow the initial random variable $\hat
X^{(n)}_0$ to have an arbitrary distribution on $\IS^{(n)}_X$. When
$X^{(n)}_0$ has distribution $\mu$ on $\IS^{(n)}_X$, we will denote the
corresponding probability and expectation operator by $P_\mu $ and
$E_\mu$, respectively. If $\mu=\delta_x$ for some \mbox{$x\in\IS^{(n)}_X$}, we
will instead write $P_x$ and $E_x$, respectively. When considering an
initial condition $x$ for $\hat X^{(n)}$, $x$ will always be in
$\IS^{(n)}_X$, although this will frequently be suppressed in the
notation. The symbols $E$ and $P$ (without any subscripts) will
correspond to the initial distribution as in Condition
\ref{conditioncr}.

\begin{lemma}\label{lemexpmomentsonfinitetime}
Suppose Conditions~\ref{conditioncr}\textup{(i)} and \textup{(ii)}, \ref
{conditionc1} and~\ref{conditionmoment} hold.
Then there exist $\delta,\rho\in(0,\infty)$ such that for every $M>0$,
%
\begin{equation}
\label{eqexponentialmomentoncompacts} \sup_{n\in\IN, x\le M} E_x
\Bigl(\sup_{0\le t\le\rho}e^{\delta
\hat X^{(n)}_t} \Bigr)=:d(\delta,\rho,M)<\infty.
\end{equation}
\end{lemma}

\begin{lemma}\label{lemsuptsupnexpmoment}
Suppose Conditions~\ref{conditioncr}\textup{(i)} and \textup{(ii)}, \ref
{conditionc1} and~\ref{conditionmoment} hold. Then there exist
$\delta,\tilde d\in(0,\infty)$ such that for every $x\in\IS^{(n)}_X$, $n\in\IN$ and
$t\ge0$,
%
\begin{equation}
\label{eqsupofexpmoment} E_x\bigl(e^{{\delta}\hat X^{(n)}_t/{2}}\bigr)\le\tilde d
e^{\delta x}.
\end{equation}
\end{lemma}
The following is immediate from Lemmas~\ref{lemEVofsquaredsuprema} and
\ref{lemsuptsupnexpmoment}.
%
\begin{corollary}\label{corEVofsquaredsupremaofshiftedhatX^n}
Suppose Conditions~\ref{conditioncr}\textup{(i)} and \textup{(ii)}, \ref
{conditionc1} and~\ref{conditionmoment} hold.
Let $\delta$ be as in Lemma~\ref{lemsuptsupnexpmoment} and
$T\in\IR_+$. Then there exists a $d(\delta,T)\in(0,\infty)$ such
that for all $x\in\IS^{(n)}_X$ and $n\in\IN$,
%
\begin{equation}
\label{eqestimateofshiftedhatX^noverfinitetime} \sup_{s\in\IR_+}
E_x \Bigl( \sup_{s\le u\le s+T}\bigl(\hat X^{(n)}_u
\bigr)^2 \Bigr)\le d(\delta,T)e^{\delta x}.\vadjust{\goodbreak}
\end{equation}
\end{corollary}

The next lemma follows by combining Lemma~\ref{lemsuptsupnexpmoment}
with arguments as in the proof of Lemma~\ref{lemEVofsquaredsuprema}.
The proof is omitted.

\begin{lemma}\label{lemEVofsquaredsupremaofshiftedprocesses}
Suppose Conditions~\ref{conditioncr}\textup{(i)} and \textup{(ii)}, \ref
{conditionc1} and~\ref{conditionmoment} hold. Let $\delta$ be as
in Lemma~\ref{lemsuptsupnexpmoment}.
Then for each $T\geq0$ there are $L_T,\tilde L_T\in(0,\infty)$ such
that for all $n\in\IN$ and $s\in\IR_+$,
%
\begin{eqnarray}
\label{eqEVofsquaredsupremaofshiftedprocesses}
&&E \Bigl(\sup_{t\le T}
\bigl(\bigl(\hat X^{(n)}_{s+t}-\hat X^{(n)}_s
\bigr)^2+\bigl(M^{(n)}_{s+t}(\phi_1)-M^{(n)}_s(
\phi_1)\bigr)^2
\nonumber
\\
&&\hspace*{42.5pt}\qquad{} +\bigl(\hat N^{(n)}_{s+t}-\hat N^{(n)}_s
\bigr)^2+\bigl(\hat\eta^{(n)}_{s+t}-\hat
\eta^{(n)}_s\bigr)^2 \bigr) \Bigr)
\\
&&\qquad\le L_T E \bigl(\hat X^{(n)}_s
\bigr)^2\le\tilde L_T \bigl(e^{\delta x^{(n)}_0} \bigr).\nonumber
\end{eqnarray}
\end{lemma}

In order to prove weak convergence results for the scaled catalyst and
reactant processes, we will need to argue that the limit processes are
continuous, which will be a consequence of the following bounds on the
jumps. The somewhat stronger estimate on the jumps of the catalyst
population in (\ref{eqPsupDeltahatX^nt}), below, will be used
in the stochastic averaging argument in the proof of Theorem \ref
{fastcatalystBPaveraging}.
Recall that for a process $\{\xi_t\}_{t\in\IR+}$ the jump at instant
$t>0$ is defined as $\Delta\xi_t:=\xi_t-\xi_{t-}$ and $\Delta\xi_0:=0$.

\begin{lemma}\label{lemPsupDeltahatX^nt}
Suppose Condition~\ref{conditioncr} holds. Fix $T,\varepsilon>0$. Then,
as \mbox{$n\to\infty$},
%
\begin{equation}
\label{eqPsupDeltahatY^nt} P \Bigl(\sup_{0\le t\le T} \bigl( \bigl|
\Delta\hat X^{(n)}_{t}\bigr|+\bigl|\Delta \hat Y^{(n)}_{t}\bigr|
\bigr)\ge\varepsilon \Bigr)\to0.
\end{equation}
If additionally Conditions~\ref{conditionc1} and \ref
{conditionmoment} hold, then, as $n\to\infty$,
%
\begin{equation}
\label{eqPsupDeltahatX^nt} \sup_{s\in\IR_+}P \Bigl(
\sup_{0\le t\le T} \bigl|\Delta\hat X^{(n)}_{s+t}\bigr|\ge\varepsilon
\Bigr)\to0.
\end{equation}
\end{lemma}

\subsection{Proofs of auxiliary results}\label{secproofsofauxresults}
In this section we prove the results stated in Section~\ref{secauxresults}.
We begin with the proofs of Lemmas \ref
{lemrepresentationofcrprocesses}--\ref{lemEVofsquaredsupremaofshiftedprocesses}.
Using these results, we will then prove Proposition \ref
{lemtighnessofcatalyst-reactant}. The proof of Lemma \ref
{lemPsupDeltahatX^nt} is given at the end.
\begin{pf*}{Proof of Lemma~\ref{lemrepresentationofcrprocesses}}
Recall that $\hat{\mathbf W}^{(n)}=(\hat X^{(n)}, \hat Y^{(n)}, \hat
Z^{(n)})$ and that for $\mathbf x=(x_1,x_2,x_3)\in\IW$, $\phi_i(\mathbf
x)=x_i$, $i=1,2,3$, and $h:=\phi_1-\phi_3$.
From (\ref{locmart}),
%
\begin{equation}
\label{eqhatZ^n} \hat Z_t^{(n)}=
\phi_3\bigl(\hat{\mathbf W}_t^{(n)}\bigr)=\hat
Z_0^{(n)}+\int_0^t
\mathcal{\hat A}^{(n)}\phi_3\bigl(\hat{\mathbf
W}_s^{(n)}\bigr)\,ds+M_t^{(n)}(
\phi_3).
\end{equation}
Using (\ref{eqgeneratorofwholeprocess}), we get
%
\begin{equation}
\label{eqhatA^nphi3} \mathcal{\hat A}^{(n)}
\phi_3\bigl(\hat{\mathbf W}_t^{(n)}\bigr)=
\lambda_1^{(n)} n \hat X_t^{(n)} \sum
_{k=0}^\infty(k-1)\mu_1^{(n)}(k)=c_1^{(n)}
\lambda_1^{(n)} \hat X_t^{(n)}.
\end{equation}
%
Next, since $\hat X_0^{(n)}=\hat Z_0^{(n)}$, we have
\[
\hat X_t^{(n)}-\hat Z_t^{(n)}=h\bigl(
\hat{\mathbf W}_t^{(n)}\bigr)=\int_0^t
\mathcal {\hat A}^{(n)}h\bigl(\hat{\mathbf W}_s^{(n)}
\bigr)\,ds+M_t^{(n)}(h)
\]
and, once more using (\ref{eqgeneratorofwholeprocess}),
\[
\mathcal{\hat A}^{(n)}h(\mathbf w)=\lambda_1^{(n)}
n \mu_1^{(n)}(0) 1_{\{x=1\}},\qquad \mathbf w=(x,y,z).
\]
Thus with $\hat\eta^{(n)}$ as in (\ref{hateta^{n}}),
we get
%
\begin{equation}
\label{eqhatX^n-hatZ^n} \hat X_t^{(n)}-
\hat Z_t^{(n)}=h\bigl(\hat{\mathbf W}^{(n)}_t
\bigr)=\hat\eta_t^{(n)}+M_t^{(n)}(h).
\end{equation}
Noting that $M^{(n)}(\phi_1)=M^{(n)}(h)+M^{(n)}(\phi_3)$ and using
(\ref{eqhatZ^n}), (\ref{eqhatA^nphi3}) and (\ref{eqhatX^n-hatZ^n}), we have
%
\begin{equation}
\label{eqhatX^nintemsofmartandeta} \hat X_t^{(n)}=
\hat X^{(n)}_0+c_1^{(n)}
\lambda_1^{(n)}\int_0^t
\hat X_s^{(n)}\,ds+M_t^{(n)}(
\phi_1)+\hat\eta_t^{(n)}.
\end{equation}
Since $\hat\eta^{(n)}$ is nondecreasing and $\int_0^\infty1_{\{
\hat X^{(n)}_s\neq1\}}\,d\hat\eta_s^{(n)}=0$, we have from the
characterization given above (\ref{eqSkorohodLipschitz}) that
\[
\hat X_t^{(n)}=\Gamma \biggl(\hat X_0^{(n)}+c_1^{(n)}
\lambda_1^{(n)}\int_0^\cdot
\hat X_s^{(n)}\,ds+M_\cdot^{(n)}(
\phi_1) \biggr) (t).
\]

Next, for the reactant population, using similar calculations as for $
\hat X^{(n)}$, we get
\begin{eqnarray*}
\hat Y^{(n)}_t&=&\hat Y^{(n)}_0+\int
_0^t \mathcal{\hat A}^{(n)}
\phi_2\bigl(\hat{\mathbf W}_s^{(n)}
\bigr)\,ds+M^{(n)}_t(\phi_2)
\\
&=&\hat Y^{(n)}_0+ c_2^{(n)}
\lambda_2^{(n)} \int_0^t
\hat X_s^{(n)}\hat Y_s^{(n)}
\,ds+M^{(n)}_t(\phi_2).
\end{eqnarray*}
%
Finally, routine calculations then show (see~\cite{Joffe1986}, Lemma
3.1.3) that (\ref{eqpqvM^nphi1}) and (\ref{eq5.6}) hold.
Details are omitted.
\end{pf*}
\begin{pf*}{Proof of Lemma~\ref{lemEVofsquaredsuprema}}
Using (\ref{eqpqvM^nphi1}) and Doob's inequality, we have
%
\begin{equation}
\label{eq5.20} E \Bigl(\sup_{t\leq T} \bigl(M^{(n)}_t(
\phi_1) \bigr)^2 \Bigr)\leq 4 \lambda_1^{(n)}
\alpha_1^{(n)}E \biggl(\int_0^{ T}
\hat X_s^{(n)}\,ds \biggr).
\end{equation}
Next, from (\ref{eq5.3}),
$
\hat X^{(n)}_t=\Gamma(\hat N^{(n)}_\cdot)(t).
$
The Lipschitz continuity of the Skorohod map implies
%
\begin{equation}
\label{eq5.21} \sup_{t\leq T}\bigl|\hat X_t^{(n)}-1\bigr|\leq2
\sup_{t\leq T}\bigl|\hat N_t^{(n)}-1\bigr|.
\end{equation}
Letting $ |\hat X^{(n)}|^2_{*,T}:=\sup_{t\le T}|\hat X^{(n)}_t|^2 $,
we now get
\[
\bigl|\hat X^{(n)}\bigr|^2_{*,T} \le2 \bigl|\hat
X^{(n)}-1\bigr|^2_{*,T}+2\le8 \bigl|\hat
N^{(n)}-1\bigr|^2_{*,T}+2\le16 \bigl|\hat
N^{(n)}\bigr|^2_{*,T}+18.
\]
Combining this with (\ref{eqhatN^ndef}) and (\ref{eq5.20}), we obtain
%
\begin{eqnarray}\label{eq5.22}
E \bigl(\bigl|\hat X^{(n)}\bigr|_{*,T}^2 \bigr)
&\le&18+ 16 E
\bigl(\bigl|\hat N^{(n)}\bigr|_{*,T}^2 \bigr)
\nonumber\\
&\le&18+48 \biggl[E\bigl(\hat X^{(n)}_0
\bigr)^2\\
&&\hspace*{39pt}{}+ \bigl(T \bigl(c^{(n)}_1
\lambda^{(n)}_1 \bigr)^2 +4
\lambda^{(n)}_1\alpha^{(n)}_1 \bigr)
\int_0^T E \bigl(\bigl|\hat X^{(n)}\bigr|_{*,s}^2
\bigr)\,ds \biggr].
\nonumber
\end{eqnarray}
Using Gronwall's inequality, we get, since $E(\hat
X^{(n)}_0)^2=(x^{(n)}_0)^2 \ge1$,
\[
E \bigl(\bigl|\hat X^{(n)}\bigr|_{*,T}^2 \bigr)\leq66
\bigl(x^{(n)}_0\bigr)^2\exp
\bigl(K^{(n)}_{1,T} \bigr),
\]
where $K^{(n)}_{1,T}:=
48T(T (c^{(n)}_1\lambda^{(n)}_1)^2+4\lambda^{(n)}_1\alpha^{(n)}_1)$.
Since $c^{(n)}_1,\lambda^{(n)}_1$ and $\alpha^{(n)}$ converge as
$n\to\infty$, we have that for some $K\in(0,\infty)$ and all $n\in
\IN$
%
\begin{equation}
\label{eq5.23} E \Bigl(\sup_{t\leq T}\bigl(\hat X_{ t}^{(n)}
\bigr)^2 \Bigr)\le66\exp\bigl(KT^2\bigr)
\bigl(x^{(n)}_0\bigr)^2.
\end{equation}
Using (\ref{eq5.23}) in (\ref{eq5.20}), (\ref
{eq5.22}) and (\ref{eqhatX^nintemsofmartandeta}), we
have the estimate in (\ref{eqEVofsquaredsuprema}) by
choosing $K$ sufficiently large.

We next establish (\ref{eqestimateonhatY^n}).
Using Doob's
inequality once more and applying (\ref{eq5.6}), we have
\begin{eqnarray*}
E \Bigl(\sup_{t\leq T} \bigl(M_{\sigma^{(n)}_{k}\wedge t}^{(n)}(
\phi_2) \bigr)^2 \Bigr) &\leq& 4 E \bigl(\bigl\langle
M^{(n)}(\phi_2)\bigr\rangle_{\sigma^{(n)}_{k}\wedge T} \bigr)
\\
&\leq& 4 \lambda_2^{(n)}\alpha_2^{(n)}E
\biggl(\int_0^{\sigma
^{(n)}_{k}\wedge T}\hat X_s^{(n)}
\hat Y_s^{(n)}\,ds \biggr).
\end{eqnarray*}
%
Thus, by (\ref{eqhatY^nb^n2}),
\begin{eqnarray*}
&&
E \bigl(\bigl|\hat Y^{(n)}\bigr|_{*,T\wedge\sigma^{(n)}_k}^2 \bigr)
\\
&&\qquad\leq3 \biggl(\bigl( y^{(n)}_0\bigr)^2+
\bigl[T \bigl(c^{(n)}_2\lambda^{(n)}_2 k
\bigr)^2+4\lambda^{(n)}_2\alpha^{(n)}_2k
\bigr]\int_0^T E \bigl(\bigl|\hat
Y^{(n)}\bigr|_{*,s\wedge\sigma^{(n)}_k}^2 \bigr)\,ds \biggr).
\end{eqnarray*}
The estimate in (\ref{eqestimateonhatY^n}) now follows by choosing
$K$ sufficiently large and applying Gronwall's inequality.
\end{pf*}
\begin{pf*}{Proof of Lemma~\ref{lemexpmomentsonfinitetime}}
First we show, using Conditions~\ref{conditioncr}(i) and (ii), \ref
{conditionc1} and~\ref{conditionmoment}, that there are $\delta_0,
d_1, d_2\in(0,\infty)$ such that for all $\delta\in[0,\delta_0]$
and~\mbox{$n\in\IN$}
%
\begin{equation}
\label{eqinStep1} -\delta d_2\le\sum_{k=0}^\infty
n^2 \bigl[e^{{(k-1)\delta
}/{n}}-1 \bigr]\mu^{(n)}_1(k)
\le-\delta d_1.
\end{equation}
Note that
\begin{eqnarray*}
&&
\sum_{k=0}^\infty n^2
\bigl[e^{{(k-1)\delta}/{n}}-1 \bigr]\mu^{(n)}_1(k)\\
&&\qquad=n^2
\sum_{k=0}^\infty \Biggl(\sum
_{l=1}^\infty\frac
{1}{l!} \biggl(
\frac{(k-1)\delta}{n} \biggr)^l \Biggr)\mu^{(n)}_1(k)
\\
&&\qquad=
n\delta \Biggl(\sum_{k=0}^\infty k
\mu^{(n)}_1(k)-1 \Biggr) + \frac
{1}{2}
\delta^2\sum_{k=0}^\infty(k-1)^2
\mu^{(n)}_1(k)
\\
&&\qquad\quad{} +n^2\sum_{k=0}^\infty \Biggl(
\sum_{l=3}^\infty\frac
{1}{l!} \biggl(
\frac{(k-1)\delta}{n} \biggr)^l \Biggr)\mu^{(n)}_1(k).
\end{eqnarray*}
Now, as $n\to\infty$,
\[
n\delta \Biggl(\sum_{k=0}^\infty k
\mu^{(n)}_1(k)-1 \Biggr) =n\delta \bigl(m^{(n)}_1-1
\bigr)=\delta c^{(n)}_1\to\delta c_1 \in(-
\infty,0)
\]
and
\[
\frac{1}{2}\delta^2\sum_{k=0}^\infty(k-1)^2
\mu^{(n)}_1(k)=\frac
{1}{2}\delta^2
\alpha^{(n)}_1\to\frac{1}{2}\delta^2
\alpha_1.
\]
Noting that $c^{(n)}_1<0$, we can choose $\delta_0>0$ sufficiently
small, and $d_1,d_2\in(0,\infty)$ suitably, such that (\ref
{eqinStep1}) holds.

For $\delta_0$ as above and $\delta\le\delta_0$, let
\[
\alpha^{(n)}_\delta:=n e^\delta\sum
_{k=1}^\infty \bigl(e^{
{(k-1)\delta}/{n}}-1 \bigr)
\frac{\mu^{(n)}_1(k) }{\mu^{(n)}_1(0)}
\]
and
\[
\beta^{(n),\delta}_t:= n^2\lambda^{(n)}_1
\int_0^t\hat X^{(n)}_s
\sum_{k=0}^\infty \bigl(
\bigl[e^{{(k-1)\delta}/{n}}-1 \bigr]\mu^{(n)}_1(k)
\bigr)1_{\{\hat X^{(n)}_s> 1\}}\,ds.
\]
Note that, by (\ref{eqinStep1}), for any $t\ge u \ge0$,
%
\begin{eqnarray}
\label{eqestimatebeta^bdelta} -\delta d_2
\lambda^{(n)}_1 \int_u^t
\hat X^{(n)}_s 1_{\{\hat
X^{(n)}_s>1\}}\,ds&\le&\beta^{(n),\delta}_t
- \beta^{(n),\delta
}_u
\nonumber\\[-8pt]\\[-8pt]
& \le&-\delta d_1\lambda^{(n)}_1 \int
_u^t \hat X^{(n)}_s
1_{\{\hat
X^{(n)}_s>1\}}\,ds.\nonumber
\end{eqnarray}
Moreover,
%
\begin{equation}
\label{eqestimatealpha^bdelta} 0\le\alpha^{(n)}_\delta
\le e^\delta\delta.
\end{equation}
The first inequality in the last display is immediate, the second
inequality can be seen as follows:
\begin{eqnarray*}
\alpha^{(n)}_\delta&=&n e^\delta\sum
_{k=1}^\infty \bigl(e^{
{(k-1)\delta}/{n}}-1 \bigr)
\frac{\mu^{(n)}_1(k)}{\mu^{(n)}_1(0)}
\\
&=&n e^\delta\sum_{k=0}^\infty
\bigl(e^{{(k-1)\delta
}/{n}}-1 \bigr)\frac{\mu^{(n)}_1(k)}{\mu^{(n)}_1(0)} -n e^\delta
\bigl(e^{-{\delta}/{n}}-1\bigr)\frac{\mu^{(n)}_1(0)}{\mu^{(n)}_1(0)}.
\end{eqnarray*}
By (\ref{eqinStep1}), the first term on the right-hand side of the
last display is smaller or equal to $0$. Thus
\[
\alpha^{(n)}_\delta\leq-n e^\delta
\bigl(e^{-{\delta}/{n}}-1\bigr)\leq e^\delta\delta.
\]

We now argue that
\[
M^{(n),\delta}_t:=\exp \bigl(\delta\hat X^{(n)}_t-
\beta^{(n),\delta
}_t \bigr)-\alpha^{(n)}_\delta
\int_0^t \exp \bigl(-\beta^{(n),\delta}_s
\bigr)\,d\hat\eta^{(n)}_s
\]
is a local martingale. Let $f(x)=e^{\delta x}$ and
\[
q(x):=\frac{\mathcal{\hat L}^{(n)} f(x)}{f(x)} 1_{\{x>1\}}.
\]
Here $\mathcal{\hat L}^{(n)}$ is the generator of $\hat X^{(n)}$, that
is, for $x \in\IS^{(n)}_X$,
%
\begin{equation}\quad
\label{eq5.27} \mathcal{\hat L}^{(n)}f(x)=\lambda_1^{(n)}
n^2 x\sum_{k=0}^\infty \biggl[f
\biggl(1 \vee \biggl( x+\frac
{k-1}{n} \biggr) \biggr)-f(x) \biggr]
\mu_1^{(n)}(k).
\end{equation}
Note that
\[
q(x)=n^2\lambda^{(n)}_1 x \sum
_{k=0}^\infty \bigl( \bigl[e^{
{(k-1)\delta}/{n}}-1 \bigr]
\mu^{(n)}_1(k) \bigr)1_{\{x>1\}}
\]
and thus
%
\begin{equation}
\label{eq5.28} \int_0^t q\bigl(\hat
X^{(n)}_s\bigr)\,ds=\beta^{(n),\delta}_t.
\end{equation}
Also,
%
\begin{equation}
\label{eq5.29} \int_0^t\mathcal{\hat
L}^{(n)}f\bigl(\hat X^{(n)}_s\bigr)1_{\{\hat X^{(n)}_s=1\}
}\,ds=
\alpha^{(n)}_\delta\hat\eta^{(n)}_t.
\end{equation}
Consider the Markov process $V^{(n)}$ defined by
\[
V^{(n)}_t:= \biggl(\hat X^{(n)}_t,\exp
\biggl(-\int_0^tq\bigl(\hat
X^{(n)}_s\bigr)\,ds \biggr) \biggr),\qquad t\ge0.\vadjust{\goodbreak}
\]
Denote by $\bar{\mathcal{L}}^{(n)}$ the generator of $V^{(n)}$. Then
the action of the generator on the function $f(x)g(y)$ with
$f(x)=e^{\delta x}$ and $g(y)=y$ is given by
\[
\bar{\mathcal{L}}^{(n)}\bigl(f(x)g(y)\bigr)=y \bigl(\mathcal{\hat
L}^{(n)}f(x)-q(x)f(x) \bigr)=y\mathcal{\hat L}^{(n)}f(x)
1_{\{x=1\}}.
\]
Using (\ref{eq5.28}) we now have that
\begin{eqnarray*}
&&
(fg)\bigl(V^{(n)}_t\bigr)-\int_0^t
\bar{\mathcal{L}}^{(n)}(fg) \bigl(V^{(n)}_s\bigr)\,ds
\\
&&\qquad=e^{\delta\hat X^{(n)}_t-\beta^{(n),\delta}_t} -\int_0^t
e^{-\beta
^{(n),\delta}_s}\mathcal{\hat L}^{(n)}f\bigl(\hat X^{(n)}_s
\bigr)1_{\{\hat
X^{(n)}_s=1\}}\,ds,\qquad t\ge0,
\end{eqnarray*}
is a local martingale.
From (\ref{eq5.29}) we now see that the last expression equals
$M^{(n),\delta}_t$, $t\ge0$, which is thus a local martingale.

We next show that for every $M>0$, $\delta\le\delta_0$ and $t\ge0$,
%
\begin{equation}
\label{eq5.30} d_3(\delta,t,M):=\sup_{x \le M}
\sup_{n\in\IN} E_x \bigl(e^{\delta
\hat X^{(n)}_t} \bigr)<\infty.
\end{equation}
%
Note that
%
\begin{eqnarray}
\label{eq5.31}
e^{\delta\hat X^{(n)}_t} &=& \biggl(e^{\delta\hat X^{(n)}_t-\beta^{(n),\delta}_t}
-\alpha^{(n)}_\delta
\int_0^t e^{-\beta^{(n),\delta}_s}\,d\hat
\eta^{(n)}_s +\alpha^{(n)}_\delta\int
_0^t e^{-\beta^{(n),\delta}_s}\,d\hat
\eta^{(n)}_s \biggr)e^{\beta^{(n),\delta}_t}\hspace*{-20pt}
\nonumber\\[-4pt]\\[-12pt]
&=& \biggl(M^{(n),\delta}_t+\alpha^{(n)}_\delta
\int_0^t e^{-\beta
^{(n),\delta}_s}\,d\hat
\eta^{(n)}_s \biggr)e^{\beta^{(n),\delta}_t}.\nonumber 
\end{eqnarray}
Applying It\^o's formula and using (\ref{eqinStep1}), (\ref
{eq5.28}) and (\ref{eq5.31}), we see that
\begin{eqnarray*}
&&
E_x \bigl(e^{\delta\hat X^{(n)}_t} \bigr)
\\
&&\qquad=e^{\delta x} + \alpha^{(n)}_\delta E_x
\int_0^t e^{\int_0^s q(\hat
X^{(n)}_u)\,du} e^{-\int_0^s q(\hat X^{(n)}_u)\,du}\,d
\hat\eta^{(n)}_s
\\
&&\qquad\quad{} + E_x \biggl(\int_0^t q
\bigl(\hat X^{(n)}_s\bigr) \biggl(M^{(n),\delta}_s
+\alpha^{(n)}_\delta\int_0^s
e^{-\beta^{(n),\delta}_u }\,d\hat\eta^{(n)}_u \biggr)e^{\int_0^s q(\hat X^{(n)}_u)\,du}
\,ds \biggr)
\\
&&\qquad= e^{\delta x} + \alpha^{(n)}_\delta E_x
\hat\eta^{(n)}_t + E_x \biggl(\int
_0^t q\bigl(\hat X^{(n)}_s
\bigr) e^{\delta\hat X^{(n)}_s} \,ds \biggr) \\
&&\qquad\le e^{\delta x} + \alpha^{(n)}_\delta
E_x\hat\eta^{(n)}_t.
\end{eqnarray*}
%
The estimate in (\ref{eq5.30}) now follows
by combining the above inequality with Lem\-ma~\ref
{lemEVofsquaredsuprema} and (\ref{eqestimatealpha^bdelta}).

Fix $M>0$, $x\le M$ and $\delta\le\frac{\delta_0}{4}$. Then, since
$\beta^{(n),\delta}_t\le0$ for all $t\ge0$, we have for $\rho>0$,
\[
E_x \Bigl(\sup_{0\le t\le\rho}
e^{\delta\hat X^{(n)}_t} \Bigr)^2\le E_x \Bigl(
\sup_{0\le t\le\rho} e^{\delta\hat
X^{(n)}_t-\beta^{(n),\delta}_t} \Bigr)^2 \le4E_x
\bigl(e^{2\delta\hat
X^{(n)}_\rho-2\beta^{(n),\delta}_\rho} \bigr), %
\]
where the last inequality follows on noting that $e^{\delta\hat
X^{(n)}_t-\beta^{(n),\delta}_t}$ is a submartingale and applying
Doob's inequality.
Now from (\ref{eqestimatebeta^bdelta}) and (\ref{eq5.30}),
\begin{eqnarray*}
E_x \bigl(e^{2\delta\hat X^{(n)}_\rho-2\beta^{(n),\delta}_\rho
} \bigr)&\le& \bigl(E_x
\bigl(e^{4\delta\hat X^{(n)}_\rho} \bigr) \bigr)^{1/2} \bigl(
E_x \bigl(e^{-4\beta^{(n),\delta}_\rho
} \bigr) \bigr)^{1/2}
\\
&\le& \bigl(d_3(4\delta,\rho,M) \bigr)^{1/2}
E_x \Bigl(\exp \Bigl(4 \delta d_2 \lambda^{(n)}_1
\rho\sup_{0\le t\le\rho}\hat X^{(n)}_t \Bigr) \Bigr).
\end{eqnarray*}
Choose $\rho< (8 d_2 \sup_{n\in\IN}\lambda^{(n)}_1
)^{-1}$. Then, by combining the above estimates, we can find a
$d_4(\delta,\rho,M)<\infty$ such that for all $x\le M$, $n\in\IN$
and $\delta\le\frac{\delta_0}{4}$
\begin{eqnarray*}
E_x \Bigl(\sup_{0\le t\le\rho}e^{\delta\hat X^{(n)}_t} \Bigr)&\le&
d_4(\delta,\rho,M)E_x \Bigl(\exp \Bigl( 4 \delta
d_2 \lambda^{(n)}_1 \rho\sup_{0\le t\le\rho}
\hat X^{(n)}_t \Bigr) \Bigr)
\\
&\le& d_4(\delta,\rho,M)E_x \biggl( \exp \biggl(
\frac{\delta}{2}\sup_{0\le
t\le\rho}\hat X^{(n)}_t \biggr)
\biggr)\\
&\le& d_4(\delta,\rho,M) \Bigl[E_x \Bigl(
\sup_{0\le t\le\rho} e^{\delta\hat
X^{(n)}_t} \Bigr) \Bigr]^{1/2}.
\end{eqnarray*}
Dividing both sides by $ [E_x  (\sup_{0\le t\le\rho}
e^{\delta\hat X^{(n)}_t} ) ]^{1/2}$ yields
\[
\Bigl[E_x \Bigl(\sup_{0\le t\le\rho} e^{\delta\hat
X^{(n)}_t} \Bigr)
\Bigr]^{1/2}\le d_4(\delta,\rho,M)
\]
for any $x\le M$ and $n\in\IN$. The result follows.
\end{pf*}
\begin{pf*}{Proof of Lemma~\ref{lemsuptsupnexpmoment}} For
$\delta\in(0,1), n\in\IN$, define
\begin{eqnarray*}
b^{(n),1}_\delta(x)&:=&\lambda^{(n)}_1
n^2 x\sum_{k=0}^\infty
\bigl(e^{\delta({k-1})/{n}}-1 \bigr)\mu^{(n)}_1(k),
\\
b^{(n),2}_\delta(x)&:=&\lambda^{(n)}_1
n^2 x\sum_{k=1}^\infty
\bigl(e^{\delta({k-1})/{n}}-1 \bigr)\mu^{(n)}_1(k)
\end{eqnarray*}
and
\[
b^{(n)}_\delta(x):=b^{(n),1}_\delta(x)
1_{\{x>1\}}+b^{(n),2}_\delta (x)1_{\{x=1\}}.
\]
From (\ref{eqinStep1}), we have, for some $\kappa\in(0,\infty)$,
\[
\sup_{n\in\IN}b^{(n),1}_\delta(x)\le-\delta
d_1x\inf_{n\in\IN} \lambda^{(n)}\le-\delta\kappa x
\le-\delta\kappa
\]
for all $\delta\le\delta_0$ [with $\delta_0$ as above (\ref
{eqinStep1})] and $x\ge1$.
Observing that with $f(x)=e^{\delta x}$, $\frac{\mathcal{\hat
L}^{(n)}f(x)}{f(x)}=b^{(n)}_\delta(x)$, where $\mathcal{\hat
L}^{(n)}$ is the generator of $\hat X^{(n)}$ defined in (\ref
{eq5.27}), we have that
%
\begin{equation}
\label{eqU^nmart} U^{(n)}_t:=e^{\delta\hat X^{(n)}_t-\int_0^t b^{(n)}_\delta(\hat
X^{(n)}_s)\,ds},\qquad t\geq0,
\end{equation}
is a local martingale.
Fix $\delta$ and $\rho$ as in the statement of Lemma \ref
{lemexpmomentsonfinitetime}. Without loss of generality, we can assume that
$\delta\le\delta_0$. Note that on the set
\[
\bigl\{\omega\dvtx  \hat X^{(n)}_s(\omega)>1\mbox{ for all } s
\in\bigl[(j-1)\rho,j\rho\bigr)\bigr\},
\]
we have
%
\begin{eqnarray}
\label{eq5.33} \qquad\delta\bigl[\hat X^{(n)}_{j\rho}-\hat
X^{(n)}_{(j-1)\rho}\bigr]&\le&\delta \bigl[\hat X^{(n)}_{j\rho}-
\hat X^{(n)}_{(j-1)\rho}\bigr]-\int_{(j-1)\rho
}^{j\rho}b^{(n)}_\delta
\bigl(\hat X^{(n)}_s\bigr)\,ds-\delta\kappa\rho
\nonumber\\[-8pt]\\[-8pt]
&\equiv& v^{(n)}_j-\delta\kappa\rho.\nonumber
\end{eqnarray}
Fix $t>0$, and let $N\in\IN$ be such that $(N-1)\rho\le t< N\rho$.
Then, similarly, on the set
\[
\bigl\{\omega\dvtx  \hat X^{(n)}_t(\omega)>1\mbox{ for all } s
\in\bigl[(N-1)\rho,t\bigr)\bigr\},
\]
$\delta[\hat X^{(n)}_{t}-\hat X^{(n)}_{(N-1)\rho}]\le
v^{(n)}_{N}(t)$, where
\[
v^{(n)}_{j}(t):=\delta\bigl[\hat X^{(n)}_{t}-
\hat X^{(n)}_{(j-1)\rho
}\bigr]-\int_{(j-1)\rho}^{t}b^{(n)}_\delta
\bigl(\hat X^{(n)}_s\bigr)\,ds.
\]
Now, for a fixed $\omega$, let $m\equiv m(\omega)$ be such that
$[(m-1)\rho, m\rho)$ is the last interval in which $\hat X^{(n)}$
visits 1 before time $N\rho$. We set $m=0$ if 1 is not visited before
time~$N\rho$. We distinguish between the cases $0<m<N$, $m=N$ and
$m=0$, where the latter corresponds to the case where 1 is not visited
before time $N\rho$.

\textit{Case} 1: $0<m<N$.

In this case
\[
\delta\hat X^{(n)}_t\le\delta\hat X^{(n)}_{m\rho}+
\sum_{j=m+1}^{N-1}\bigl(v^{(n)}_j-
\delta\kappa\rho\bigr)+v^{(n)}_N(t).
\]
For $j\in\IN$, let
\[
\gamma^{(n)}_j:=\inf\bigl\{t\ge(j-1)\rho| \hat
X^{(n)}_t=1\bigr\}\wedge j\rho
\]
and
%
\begin{equation}
\label{eqtheta^nj} \theta^{(n)}_j:=
\sup_{0\le t\le\rho}\bigl[\hat X^{(n)}_{(t+\gamma
^{(n)}_j)\wedge j\rho}-\hat
X^{(n)}_{\gamma^{(n)}_j}\bigr].
\end{equation}
Then
$
\delta\hat X^{(n)}_{m\rho}\le\delta\theta^{(n)}_m+\delta.
$
Combining the above estimates, we have
%
\begin{equation}
\label{eqcase1m<N} \delta\hat X^{(n)}_t\le\delta
\theta^{(n)}_m+\delta+\sum_{j=m+1}^{N-1}
\bigl(v^{(n)}_j-\delta\kappa\rho\bigr)+v^{(n)}_N(t).
\end{equation}
Thus, in this case
\[
\delta\hat X^{(n)}_t\le\delta\hat X^{(n)}_0
+\max_{0\le l \le
N} \Biggl\{\sum_{j=l+1}^{N-1}
\bigl(v^{(n)}_j-\delta\kappa\rho\bigr)+\delta
\theta^{(n)}_l \Biggr\}+v^{(n)}_N(t),
\]
where by convention $\sum_{j=l+1}^{N-1}(v^{(n)}_j-\delta\kappa\rho
)=0$ for $l=N-1, N$ and $\theta^{(n)}_0:=0$.

\textit{Case} 2: $m=0$.

In this case, 1 is not visited before time $N\rho$ and thus
\begin{eqnarray*}
\delta\hat X^{(n)}_t&\le&\delta\hat X^{(n)}_0+
\sum_{j=1}^{N-1}\bigl(v^{(n)}_j-
\delta\kappa\rho\bigr)+v^{(n)}_N(t)
\\
&\le&\delta\hat X^{(n)}_0 +\max_{0\le l \le N} \Biggl\{
\sum_{j=l+1}^{N-1}\bigl(v^{(n)}_j-
\delta\kappa\rho\bigr)+\delta\theta^{(n)}_l \Biggr
\}+v^{(n)}_N(t).
\end{eqnarray*}

\textit{Case} 3: $m=N$.

Suppose first that there is an $s\in[(N-1)\rho,t]$ such that $\hat
X^{(n)}_s=1$. It then follows that
\[
\delta\hat X^{(n)}_t\le\delta\theta^{(n)}_N+
\delta.
\]
Now suppose that there is no such $s\in[(N-1)\rho,t]$. Define $m'\in
\{1,2,\ldots,\break N-1\}$ to be such that $[(m'-1)\rho,m'\rho)$ is the
last interval in which $\hat X^{(n)}$ visits 1 before $(N-1)\rho$.
Once again we set $m'=0$ if there is no such interval.

If $m'=0$, we get exactly as in case 2 that
\[
\delta\hat X^{(n)}_t\le\delta\hat X^{(n)}_0
+\max_{0\le l \le
N} \Biggl\{\sum_{j=l+1}^{N-1}
\bigl(v^{(n)}_j-\delta\kappa\rho\bigr)+\delta
\theta^{(n)}_l \Biggr\}+v^{(n)}_N(t).
\]

If $1\le m'\le N-1$, then
\begin{eqnarray*}
\delta\hat X^{(n)}_t&\le&\delta\theta^{(n)}_{m'}+
\delta+\sum_{j=m'+1}^{N-1}\bigl(v^{(n)}_j-
\delta\kappa\rho\bigr)+v^{(n)}_N(t)
\\
&\le&\delta\hat X^{(n)}_0 +\max_{0\le l \le N} \Biggl\{
\sum_{j=l+1}^{N-1}\bigl(v^{(n)}_j-
\delta\kappa\rho\bigr)+\delta\theta^{(n)}_l \Biggr
\}+v^{(n)}_N(t).
\end{eqnarray*}

Combining the three cases, we have
\[
\delta\hat X^{(n)}_t \le\max\Biggl\{ \delta\hat
X^{(n)}_0 +\max_{ l \le N}\Biggl\{\sum
_{j=l+1}^{N-1}\bigl(v^{(n)}_j-
\delta\kappa\rho\bigr)+\delta\theta^{(n)}_l\Biggr\}
+v^{(n)}_N(t), \delta+\delta\theta^{(n)}_N
\Biggr\}.
\]
Thus, for any $M_0>0$,
\begin{eqnarray*}
&&
P_x\bigl(\delta\hat X^{(n)}_t\ge
M_0\bigr)
\\
&&\qquad\le\sum_{l=0}^{N-1} P_x
\Biggl(v^{(n)}_N(t)+\sum_{j=l+1}^{N-1}
\bigl(v^{(n)}_j-\delta\kappa\rho\bigr)+\delta
\theta^{(n)}_l+\delta+\delta \hat X^{(n)}_0
\ge M_0 \Biggr)
\\
&&\qquad\quad{} +P_x\bigl(\delta\theta^{(n)}_N+\delta\ge
M_0\bigr)
\\
&&\qquad\le e^{\delta(1+x)-M_0}
\\
&&\qquad\quad{} \times\Biggl(\sum_{l=0}^{N-1} \Biggl[
E_x \Biggl(\exp \Biggl[\delta\theta^{(n)}_l+
\sum_{j=l+1}^{N-1}v^{(n)}_j+v^{(n)}_N(t)
\Biggr] \Biggr) e^{-\delta\kappa\rho(N-l-1)} \Biggr]\\
&&\qquad\quad\hspace*{227.5pt}{}+ E_x\bigl(e^{\delta\theta
^{(n)}_N}
\bigr) \Biggr).
\end{eqnarray*}
Recalling $U^{(n)}$ from (\ref{eqU^nmart}) and using its martingale
property, we get
%
\begin{eqnarray}
\label{eq5.36}\quad P_x\bigl(\delta\hat X^{(n)}_t
\ge M_0\bigr)&\le& e^{-M_0}e^{\delta(1+x)}
\Biggl(E_x \bigl(e^{\delta\theta^{(n)}_N} \bigr)+ \sum
_{l=0}^{N-1} e^{-\delta\kappa\rho(N-l-1)} E_x
\bigl(e^{\delta
\theta^{(n)}_l} \bigr) \Biggr)
\nonumber\\[-8pt]\\[-8pt]
&\le& e^{-M_0}e^{(1+x)\delta} \,d(\delta,\rho,1) \biggl(1+
\frac
{1}{1-e^{-\delta\kappa\rho}} \biggr),\nonumber
\end{eqnarray}
where the last inequality follows from Lemma \ref
{lemexpmomentsonfinitetime} and the observation that
%
\begin{equation}
\label{eqexptailsofhatX^n} \sup_{n\in\IN}E_x
\bigl(e^{\delta\theta^{(n)}_l}\bigr)\le\sup_{n\in\IN} E_{1} \Bigl(
\sup_{0\le t\le\rho} e^{\delta\hat X^{(n)}_t} \Bigr)\le d(\delta,\rho,1)<\infty,
\end{equation}
where $\theta^{(n)}_l$ is as in (\ref{eqtheta^nj}).
Finally, from (\ref{eq5.36}), we get that
for all $t\ge0$ and $n\in\IN$
%
\begin{eqnarray}
\label{eq5.38}\quad E_x \bigl(e^{{\delta}\hat X^{(n)}_t/{2}} \bigr)&=&\int
_0^\infty P_x \bigl(\delta\hat
X^{(n)}_t>2\ln( y) \bigr)\,dy
\nonumber
\\
&\le& 1+ e^{(1+x)\delta} \,d(\delta,\rho,1) \biggl(1+\frac{1}{1-e^{-\delta\kappa\rho}} \biggr)
\int_1^\infty e^{-2\ln(y)}\,dy\\
&\le&\tilde d
e^{\delta x},\nonumber
\end{eqnarray}
where $\tilde d= 1+e^{\delta}
d(\delta,\rho,1)  (1+\frac{1}{1-e^{-\delta\kappa\rho}} )$.
The result follows.
\end{pf*}
\begin{pf*}{Proof of Proposition~\ref{lemtighnessofcatalyst-reactant}}
We will first consider the second part of the proposition.
We begin by showing that $\{\hat N^{(n)}_{s+\cdot}-\hat N^{(n)}_s\}_{s,n}$ is tight.
For that, in view of (\ref{eqEVofsquaredsupremaofshiftedprocesses}), it
suffices to show that the following condition
(Aldous--Kurtz criterion) holds: for each $M>0,\varepsilon>0$ and $\gamma
>0$ there are $\delta_0>0$ and $n_0$ such that for all stopping times
$\{\tau_n\}_{n\in\IN}$ with $\tau_n\leq M$, we have
%
\begin{equation}
\label{eq5.39} \sup_{s\in\IR_+,n\geq n_0}\sup_{\theta\leq\delta_0} P\bigl(\bigl|\hat
N^{(n)}_{s+\tau_n+\theta}-\hat N^{(n)}_{s+\tau_n}\bigr|\geq\gamma
\bigr)\leq \varepsilon.
\end{equation}

Let $M, \varepsilon,\gamma\in(0,\infty)$ be given.
Note that
\begin{eqnarray*}
&& P\bigl(\bigl|\hat N^{(n)}_{s+\tau_n+\theta}-\hat N^{(n)}_{s+\tau_n}\bigr|
\geq \gamma\bigr)
\\
&&\qquad\leq P \biggl(\biggl|c_1^{(n)} \lambda_1^{(n)}
\int_{s+\tau
_n}^{s+\tau_n+\theta}\hat X^{(n)}_u
\,du\biggr|\geq\frac{\gamma}{2} \biggr)\\
&&\qquad\quad{}+P \biggl(\bigl|M^{(n)}_{s+\tau_n+\theta}(
\phi_1)-M^{(n)}_{s+\tau_n}(\phi_1)\bigr|\geq
\frac{\gamma}{2} \biggr).
\end{eqnarray*}
By (\ref{eqestimateofshiftedhatX^noverfinitetime}) we have,
for $\delta_0$ sufficiently small,
\[
\sup_{s\in\IR_+,n\in\IN}\sup_{\theta\leq\delta_0} P \biggl(\biggl|c_1^{(n)}
\lambda_1^{(n)}\int_{s+\tau_n}^{s+\tau_n+\theta}
\hat X^{(n)}_u \,du\biggr|\geq\frac{\gamma}{2} \biggr)<
\frac{\varepsilon}{2}.
\]
It remains to prove that, for some $\delta_0>0$,
%
\begin{equation}
\label{eq5.40} \sup_{s\in\IR_+,n\in\IN}\sup_{\theta\leq\delta_0} P \biggl(\bigl|M^{(n)}_{s+\tau_n+\theta}(
\phi_1)-M^{(n)}_{s+\tau_n}(\phi_1)\bigr|\geq
\frac{\gamma}{2} \biggr) < \frac{\varepsilon}{2}.
\end{equation}
Using the martingale property of $M^{(n)}(\phi_1)$,
\begin{eqnarray*}
&&
P \biggl(\bigl|M^{(n)}_{s+\tau_n+\theta}(\phi_1)-M^{(n)}_{s+\tau_n}(
\phi_1)\bigr|\geq\frac{\gamma}{2} \biggr)\\
&&\qquad\leq\frac{E(|M^{(n)}_{s+\tau
_n+\theta}(\phi_1)-M^{(n)}_{s+\tau_n}(\phi_1)|^2)}{(\gamma/2)^2}
\\
&&\qquad=\frac{E((M^{(n)}_{s+\tau_n+\theta}(\phi_1))^2)-E((M^{(n)}_{s+\tau
_n}(\phi_1))^2)}{(\gamma/2)^2}
\\
&&\qquad=\frac{E\langle M^{(n)}(\phi_1)\rangle_{s+\tau_n+\theta}-E\langle
M^{(n)}(\phi_1)\rangle_{s+\tau_n}}{(\gamma/2)^2},
\end{eqnarray*}
and, using (\ref{eqpqvM^nphi1}),
\[
E\bigl\langle M^{(n)}(\phi_1)\bigr\rangle_{s+\tau_n+\theta}-E
\bigl\langle M^{(n)}(\phi_1)\bigr\rangle_{s+\tau_n}\le
E \biggl(\lambda_1^{(n)}\alpha_1^{(n)}
\int_{s+\tau_n}^{s+\tau_n+\theta}\hat X^{(n)}_u
\,du \biggr).
\]
Now, using (\ref{eqestimateofshiftedhatX^noverfinitetime})
once more, we can choose $\delta_0>0$ such that (\ref{eq5.40}) holds.
This proves tightness of $\{\hat N^{(n)}_{s+\cdot}-\hat N^{(n)}_s\}_{s,n}$ and,
using the continuity property of the Skorohod map (from $D_1(\IR_+\dvtx \IR
)$ to $D(\IR_+\dvtx [1,\infty))$), that of $\{\hat X^{(n)}_{s+\cdot}-\hat
X^{(n)}_s\}_{s,n}$ and $\{\hat\eta^{(n)}_{s+\cdot}-\hat\eta^{(n)}_{s}\}_{s,n}$. Tightness of $\{\hat X^{(n)}_{s+\cdot}\}_{s,n}$
now follows by using the uniform estimate in Lemma~\ref{lemsuptsupnexpmoment}.\vadjust{\goodbreak}

Now we consider the first part of the proposition. Tightness of $(\hat
X^{(n)},\hat\eta^{(n)})$ follows as before.
We now consider $\hat Y^{(n)}$.
Fix $\varepsilon>0$. Using (\ref{eqestimateonhatY^n}), we get, for
$\tilde K\in(0,\infty)$,
\begin{eqnarray*}
P \Bigl(\sup_{t\leq T}\bigl(\hat Y^{(n)}_t\bigr)>
\tilde K \Bigr)&\leq& P \Bigl(\sup_{t\leq T}\bigl(\hat Y^{(n)}_{\sigma^{(n)}_k\wedge t}
\bigr)>\tilde K \mbox { and }\sigma^{(n)}_k>T \Bigr)+P \bigl(
\sigma^{(n)}_k\leq T \bigr)
\\[-2pt]
&\leq&\frac{E (\sup_{t\leq T}(\hat Y^{(n)}_{\sigma
^{(n)}_k\wedge t})^2 )}{\tilde K^2}+P \Bigl(\sup_{t\leq T}\bigl(\hat
X^{(n)}_t\bigr)\geq k \Bigr)
\\[-2pt]
&\leq&\frac{\exp(KT^2 k^2) (y^{(n)}_0)^2}{\tilde K^2}+\frac{E
(\sup_{t\leq T}(\hat X^{(n)}_t)^2 )}{k^2}.
\end{eqnarray*}
Using (\ref{eqEVofsquaredsuprema}), we can choose $k$ such that
%
\begin{equation}
\label{eqsupE|X^n|k} \sup_{n\in\IN}\frac{E (\sup_{t\leq T}(\hat X^{(n)}_t)^2
)}{k^2}<
\frac{\varepsilon}{2}.
\end{equation}
Now choose $\tilde K $ such that
\[
\sup_{n\in\IN}\frac{\exp(KT^2 k^2)(y^{(n)}_0)^2}{\tilde K^2}<\frac{\varepsilon}{2}.
\]
The last two displays imply
$\sup_{n\in\IN} P(\sup_{t\leq T}(\hat Y^{(n)}_t)>\tilde K
)<\varepsilon$, and since $\varepsilon>0$ is arbitrary,
the tightness of the random variables $\{\hat Y^{(n)}_t\}_{n\in\IN}$,
for each $t\geq0$, follows. To establish the tightness of the
processes $\{\hat Y^{n)}\}_{n\in\IN}$, it now suffices to show that
for each $M>0,\varepsilon>0$ and $\gamma>0$ there are $\delta_0>0$ and
$n_0$ such that for all stopping times $\{\tau_n\}_{n\in\IN}$ with
$\tau_n\leq M$, we have
%
\begin{equation}
\label{eq5.42} \sup_{n\geq n_0}\sup_{\theta\leq\delta_0} P \bigl(\bigl|\hat
Y^{(n)}_{\tau_n+\theta}-\hat Y^{(n)}_{\tau_n}\bigr|\geq\gamma
\bigr)\leq\varepsilon.
\end{equation}
Fix $M,\varepsilon,\gamma\in(0,\infty)$. Then, for any $\theta\in(0,1)$,
\begin{eqnarray*}
&&
P\bigl(\bigl|\hat Y^{(n)}_{\tau_n+\theta}-\hat Y^{(n)}_{\tau_n}\bigr|
\geq \gamma \bigr)
\\[-2pt]
&&\qquad \leq P \bigl(\bigl|\hat Y^{(n)}_{(\tau_n+\theta)\wedge\sigma
^{(n)}_k}-\hat Y^{(n)}_{\tau_n\wedge\sigma^{(n)}_k}\bigr|
\geq\gamma \bigr)+P \bigl(\sigma^{(n)}_k\leq M+1 \bigr).
\end{eqnarray*}
Taking $T=M+1$ and $k$ as in (\ref{eqsupE|X^n|k}), we have
$P(\sigma^{(n)}_k<M+1)<\varepsilon/2$ for all $n\in\IN$.
For the first term on the right-hand side of the last display, we get,
using (\ref{eqhatY^nb^n2}) and that $\sup_{t\leq T\wedge\sigma
^{(n)}_k}\hat X^{(n)}_t\leq k$,
%
\begin{eqnarray}
\label{eq5.43} &&
P \bigl(\bigl|\hat Y^{(n)}_{(\tau_n+\theta)\wedge\sigma^{(n)}_k}-\hat
Y^{(n)}_{\tau_n\wedge\sigma^{(n)}_k}\bigr|\geq\gamma \bigr)
\nonumber
\\[-2pt]
&&\qquad \leq P \biggl(\biggl|c_2^{(n)}\lambda^{(n)}_2
\int_{\tau
_n\wedge\sigma^{(n)}_k}^{(\tau_n+\theta)\wedge\sigma^{(n)}_k}\hat Y^{(n)}_s
\,ds\biggr|\geq\frac{\gamma}{2k} \biggr)
\\[-2pt]
&&\qquad\quad{} +P \biggl(\bigl|M^{(n)}_{(\tau_n+\theta)\wedge
\sigma^{(n)}_k}(\phi_2)-M^{(n)}_{\tau_n\wedge\sigma^{(n)}_k}(
\phi_2)\bigr|\geq\frac{\gamma}{2} \biggr).\nonumber
\end{eqnarray}
The first term on the right-hand side can be bounded as follows:
\begin{eqnarray*}
&&
P \biggl( \biggl|c_2^{(n)}\lambda^{(n)}_2
\int_{\tau_n\wedge\sigma
^{(n)}_k}^{(\tau_n+\theta)\wedge\sigma^{(n)}_k}\hat Y^{(n)}_s
\,ds \biggr|\geq\frac{\gamma}{2k} \biggr)
\\
&&\qquad \leq \biggl(\frac{2k c_2^{(n)}\lambda^{(n)}_2}{\gamma} \biggr)^2 E \biggl( \biggl(\int
_{\tau_n\wedge\sigma^{(n)}_k}^{(\tau_n+\theta
)\wedge\sigma^{(n)}_k}\hat Y^{(n)}_s \,ds
\biggr)^2 \biggr)
\\
&&\qquad\leq\theta \biggl(\frac{2k c_2^{(n)}\lambda^{(n)}_2}{\gamma} \biggr)^2\exp \bigl(
K(M+1)^2 k^2 \bigr) \bigl(y^{(n)}_0
\bigr)^2,
\end{eqnarray*}
where $K$ is the constant from (\ref{eqestimateonhatY^n}).
Thus, for $\delta_0$ sufficiently small, we get
\[
\sup_{n\in\IN}\sup_{\theta\leq\delta_0} P \biggl(\biggl\llvert
c_2^{(n)}\lambda^{(n)}_2\int
_{\tau_n\wedge\sigma^{(n)}_k}^{(\tau
_n+\theta)\wedge\sigma^{(n)}_k}\hat Y^{(n)}_s \,ds
\biggr\rrvert \geq\frac
{\gamma}{2k} \biggr)<\varepsilon/4.
\]
The second term on the right-hand side of (\ref
{eq5.43}) can be bounded as follows:
\begin{eqnarray*}
&&
P \biggl( \bigl|M^{(n)}_{(\tau_n+\theta)\wedge\sigma^{(n)}_k}(\phi_2)-M^{(n)}_{\tau_n\wedge\sigma^{(n)}_k}(
\phi_2) \bigr|\geq\frac
{\gamma}{2} \biggr)
\\
&&\qquad\leq\frac{E ( \langle M^{(n)}(\phi_2) \rangle_{(\tau_n+\theta)\wedge\sigma^{(n)}_k}- \langle M^{(n)}(\phi_2) \rangle_{\tau_n\wedge\sigma^{(n)}_k} )}{(\gamma
/2)^2}
\\
&&\qquad\leq\frac{4}{\gamma^2} \lambda^{(n)}_2
\alpha_2^{(n)} E \biggl( \int_{\tau_n\wedge\sigma^{(n)}_k}^{(\tau_n+\theta)\wedge\sigma
^{(n)}_k}
\hat X^{(n)}_s\hat Y^{(n)}_s \,ds
\biggr)
\\
&&\qquad\leq\frac{4}{\gamma^2} \lambda^{(n)}_2
\alpha_2^{(n)} k \theta E \Bigl(\sup_{0\le s\le M+1} \hat
Y^{(n)}_{s\wedge\sigma
^{(n)}_k} \Bigr).
\end{eqnarray*}
Using (\ref{eqestimateonhatY^n}) once more, we have that, for
$\delta_0$ sufficiently small, the second term in (\ref
{eq5.43}) is bounded by $\frac{\varepsilon}{4}$.
Combining the above estimates, we now see that (\ref
{eq5.42}) holds, and thus tightness of $\{\hat Y^{(n)}\}_{n\in\IN}$ follows.
\end{pf*}
\begin{pf*}{Proof of Lemma~\ref{lemPsupDeltahatX^nt}}
Consider (\ref{eqPsupDeltahatX^nt}).
Let $N^{(n)}_{s,T}$
be the number of deaths of particles of the (unscaled) process
$X^{(n)}$ in the time interval $[s,s+T]$.
Fix $\varepsilon,\delta>0$. Then
\begin{eqnarray*}
&&
P \Bigl(\sup_{0\le t\le T} \bigl|\Delta\hat X^{(n)}_{s+t}\bigr|\ge
\varepsilon \Bigr)
\\
&&\qquad\le P \Bigl(\sup_{0\le t\le T} \bigl|\Delta\hat X^{(n)}_{s+t}\bigr|
\ge \varepsilon;\sup_{0\le t\le T} X^{(n)}_{s+t}\le nL \Bigr)+ P
\Bigl(\sup_{0\le t\le T} X^{(n)}_{s+t} > n L \Bigr).
\end{eqnarray*}
By Corollary~\ref{corEVofsquaredsupremaofshiftedhatX^n}, we
can choose
$L\in(0,\infty)$ such that
\[
P \Bigl(\sup_{0\le t\le T} X^{(n)}_{s+t} > n L \Bigr)<
\frac{\delta}{3}
\]
for $s\in\IR_+$ and $n\in\IN$. Next, consider
\begin{eqnarray*}
&&
P \Bigl(\sup_{0\le t\le T} \bigl|\Delta\hat X^{(n)}_{s+t}\bigr|\ge
\varepsilon; \sup_{0\le t\le T} X^{(n)}_{s+t}\le nL \Bigr)
\\
&&\qquad\le P \Bigl(\sup_{ t\le T} \bigl|\Delta\hat X^{(n)}_{s+t}\bigr|
\ge\varepsilon; N^{(n)}_{s,T}< nCL \Bigr)+ P \Bigl(
\sup_{ t\le T} X^{(n)}_{s+t}\le nL; N^{(n)}_{s,T}
\ge nCL \Bigr).
\end{eqnarray*}
Note that on the set $\{\sup_{0\le t\le T} X^{(n)}_{s+t}\le nL \}$ the
branching rates of $X^{(n)}$ are bounded during the time interval
$[s,s+T]$, uniformly in $s$ and $n$, and thus we can choose a $C\in
(0,\infty)$ such that for $s\in\IR_+$ and $n\in\IN$
\[
P \Bigl(\sup_{0\le t\le T} X^{(n)}_{s+t}\le
nL;N^{(n)}_{s,T} \ge nCL \Bigr)<\frac{\delta}{3}.
\]
Finally, let, for $n\in\IN$, $\{\xi^{(n)}_i\}_{i\in\IN}$ be
i.i.d. random variables distributed as $\mu^{(n)}_1$. Then, since the
variance of the offspring distribution converges, we have for $n_0$
sufficiently large and all $n\ge n_0$,
\begin{eqnarray*}
P \Bigl(\sup_{0\le t\le T} \bigl|\Delta\hat X^{(n)}_{s+t}\bigr|\ge
\varepsilon; N^{(n)}_{s,T}< nCL \Bigr)&\le& P \biggl(
\max_{1\le i<
nCL}\frac{|\xi^{(n)}_i-1|}{n} \ge\varepsilon \biggr)
\\
&\le&\sum_{i=1}^{nCL-1}P \bigl(\bigl|
\xi^{(n)}_i-1\bigr|\ge n\varepsilon \bigr)\\
&\le&\sum
_{i=1}^{nCL-1} \frac{E (|\xi^{(n)}_i-1|^2
)}{(n\varepsilon)^2}<\frac{\delta}{3}.
\end{eqnarray*}
Combining the above estimates, (\ref{eqPsupDeltahatX^nt}) follows.
The limit in (\ref{eqPsupDeltahatY^nt}) can be established
similarly, using Lemma~\ref{lemEVofsquaredsuprema} instead of
Corollary~\ref{corEVofsquaredsupremaofshiftedhatX^n}; the
proof is therefore omitted.
\end{pf*}

\section{\texorpdfstring{Proof of Theorem \protect\ref{thmcrdiffusion}}
{Proof of Theorem 2.1}}
\label{sec6}

The following martingale characterization result will be useful in the
proof of Theorem~\ref{thmcrdiffusion}. The proof is standard and is
omitted; see~\cite{Stroock1979},~\cite{Kurtz1990},~\cite{Kurtz1991}
and Theorem 5.3
of~\cite{Budhiraja2003}.

For $\phi\in C^\infty_c([1,\infty)\times\IR_+)$, let
\begin{eqnarray*}
\mathcal{L} \phi(x,y)&:=&c_1\lambda_1 x
\,\frac{\partial}{\partial
x}\phi(x,y)+\frac{1}{2}\alpha_1
\lambda_1 x \,\frac{\partial^2}{\partial x^2}\phi(x,y)
\\
&&{} +c_2 \lambda_2 x y \,\frac{\partial}{\partial y}\phi (x,y)+
\frac{1}{2}\alpha_2\lambda_2 x y \,\frac{\partial^2}{\partial
y^2}
\phi(x,y).
\end{eqnarray*}
Let $\tilde\Omega:=D(\IR_+\dvtx [1,\infty)\times\IR_+^2)$ and $\tilde
{\mathcal{F}}$ be the corresponding Borel $\sigma$-field (with
respect to the Skorohod topology). Denote by
$\{\mathcal{F}_t\}_{t\in\IR_+}$ the canonical filtration on
$(\tilde{\Omega},\tilde{\mathcal{F}})$,
that is, $\mathcal{F}_t=\sigma(\pi_s|s\leq t)$, where $\pi_s(\tilde
\omega)=\tilde\omega_s=\tilde\omega(s)$ for $\tilde\omega\in
\tilde\Omega$.
Finally, let $\pi^{(i)}$, $i=1,2,3$, be the coordinate processes, that
is, $(\pi^{(1)}(\tilde\omega),\pi^{(2)}(\tilde\omega),\pi^{(3)}(\tilde\omega))=\pi(\tilde\omega)$.
%
\begin{theorem}\label{thmmartingalecharacterization}
Let $\tilde P$ be a probability measure on $(\tilde\Omega,\tilde
{\mathcal{F}})$ under which the following hold a.s.:
\begin{longlist}[(iii)]
\item[(i)] $\pi^{(3)}$ is a nondecreasing, continuous process, and
$\pi^{(3)}_0=0$;
\item[(ii)] $(\pi^{(1)},\pi^{(2)})$ is an $([1,\infty)\times\IR_+)$ valued continuous process;
\item[(iii)] $\int_0^\infty1_{(1,\infty)}(\pi^{(1)}_s)\,d\pi^{(3)}_s=0$;
\item[(iv)] for all $\phi\in C^\infty_c([1,\infty)\times\IR_+)$
\[
\phi\bigl(\pi^{(1)}_t,\pi^{(2)}_t
\bigr)-\int_0^t\mathcal{L}\phi\bigl(
\pi^{(1)}_s,\pi^{(2)}_s\bigr)\,ds-\int
_0^t\frac{\partial\phi}{\partial
x}\bigl(1,
\pi^{(2)}_s\bigr)\,d\pi^{(3)}_s
\]
is an $\{\mathcal{F}_t\}$ martingale;
\item[(v)] $\tilde P\circ(\pi^{(1)}_0,\pi^{(2)}_0)^{-1}=\bar P\circ
(X_0,Y_0)^{-1}$,where $X,Y$ and $\bar P$ are as in Proposition \ref
{propuniqunessofX,Y}.
\end{longlist}
Then $\tilde P\circ(\pi^{(1)},\pi^{(2)})^{-1}=\bar P\circ(X,Y)^{-1}$.
\end{theorem}
\begin{pf*}{Proof of Theorem~\ref{thmcrdiffusion}}
Recall that for $\phi\in C^\infty_c([1,\infty)\times\IR_+)$, we have
%
\begin{equation}
\label{eq6.1}\quad \phi\bigl(\hat X_t^{(n)},\hat
Y_t^{(n)}\bigr)=\phi\bigl(\hat X_0^{(n)},
\hat Y_0^{(n)}\bigr)+\int_0^t
\mathcal{\hat A}^{(n)}\phi\bigl(\hat X_s^{(n)},
\hat Y_s^{(n)}\bigr)\,ds+M_t^{(n)}(\phi),
\end{equation}
where $M_t^{(n)}(\phi)$ is a martingale, and $\hat{\mathcal{
A}}^{(n)} $ is as defined in (\ref{eqgeneratorofwholeprocess}).
Also note that $\hat{\mathcal{ A}}^{(n)} $ can be rewritten as
\[
\mathcal{\hat A}^{(n)}\phi(x,y) = \mathcal{L}^{(n)}\phi(x,y)+
\mathcal{D}^{(n)}\phi(y)n\lambda_1^{(n)}
\mu_1^{(n)}(0)1_{\{x=1\}},
\]
where
\begin{eqnarray*}
\mathcal{L}^{(n)}\phi(x,y)&:=&\lambda_1^{(n)}
n^2 x\sum_{k=0}^\infty \biggl[
\phi \biggl(x+\frac{k-1}{n},y \biggr)-\phi(x,y) \biggr]\mu_1^{(n)}(k)
\\
&&{} +\lambda_2^{(n)} n^2 x y \sum
_{k=0}^\infty \biggl[\phi \biggl(x,y+\frac{k-1}{n}
\biggr)-\phi(x,y) \biggr]\mu_2^{(n)}(k)
\end{eqnarray*}
and
\[
\mathcal{D}^{(n)}\phi(y):=n \biggl(\phi(1,y)-\phi \biggl(1-
\frac
{1}{n},y \biggr) \biggr).
\]
Thus, using (\ref{hateta^{n}}), (\ref{eq6.1}) can be rewritten as
\begin{eqnarray*}
\phi\bigl(\hat X_t^{(n)}, \hat Y_t^{(n)}
\bigr) &=& \phi\bigl(\hat X_0^{(n)},\hat Y_0^{(n)}
\bigr)+\int_0^t\mathcal{L}^{(n)}\phi
\bigl(\hat X_s^{(n)},\hat Y_s^{(n)}
\bigr)\,ds
\\
&&{} +\int_0^t\mathcal{D}^{(n)}\phi
\bigl(\hat Y_{s}^{(n)}\bigr)\,d\hat\eta^{(n)}_s+M_t^{(n)}(
\phi).
\end{eqnarray*}
Recall the path space $(\tilde\Omega,\tilde{\mathcal{F}})$
introduced above Theorem~\ref{thmmartingalecharacterization}.
Denote by $\tilde{P}^{(n)}$ the measure induced by $(\hat X^{(n)},\hat
Y^{(n)}, \hat\eta^{(n)})$ on
$(\tilde{\Omega},\tilde{\mathcal{F}})$ and by $\tilde{E}^{(n)}$
the corresponding expectation.

From Proposition~\ref{lemtighnessofcatalyst-reactant}, $\tilde
{P}^{(n)}$ is tight. Let $\tilde P$ be a limit point of $\{\tilde
P^{(n)}\}$ along some subsequence $\{n_k\}$. In order to complete the
proof, it suffices to show that under $\tilde P$ properties (i)--(v) in
Theorem~\ref{thmmartingalecharacterization} hold almost surely.
Property~(i) is immediate from the fact that $\hat\eta^{(n)}$ is
nondecreasing and continuous with initial value 0 for each~$n$. Also,
property (v) is immediate from the fact that $(\hat X^{(n)}_0,\hat
Y^{(n)}_0)=(1,1)$, a.s., for each $n$.
Next, consider property (ii).
The continuity of $\pi^{(1)}$ and $\pi^{(2)}$ follows by (\ref
{eqPsupDeltahatY^nt}); see~\cite{Jacod1987}, Proposition
VI.3.26, page 315.

To see (iii), consider, for $\delta>0$,
continuous bounded test functions $f_\delta\dvtx\break [1,\infty)\longrightarrow
\IR_+$ such that
%
\begin{equation}
\label{eqfdelta} f_\delta(x) = %
\cases{1, &\quad if $x\geq1+2
\delta$,
\cr
0, &\quad if $x\leq1+\delta$. } %
\end{equation}
Note that, for each $n\in\IN$,
$\int_0^\infty f_\delta(\hat X^{(n)}_s) \,d\hat\eta^{(n)}_s=0$ and
thus, for each $\delta>0$,
\[
0=\lim_{k\to\infty}\tilde E^{(n_k)} \biggl(\int_0^\infty
f_\delta \bigl(\pi^{(1)}_s\bigr)\,d
\pi^{(3)}_s\wedge1 \biggr)=\tilde E \biggl(\int
_0^\infty f_\delta\bigl(
\pi^{(1)}_s\bigr)\,d\pi^{(3)}_s\wedge1
\biggr).
\]
Consequently, for each $\delta>0$,
$
\int_0^\infty1_{[1+2\delta,\infty)}(\pi^{(1)}_s)\,d\pi^{(3)}_s=0$,
almost surely w.r.t.~$\tilde P$. The property in (iii) now follows on
sending $\delta\to0$.

Finally, we consider part (iv). It suffices to show that for every
$0\le s\le t<\infty$
\begin{eqnarray*}
&&\tilde E \biggl( \psi(\cdot) \biggl(\phi\bigl(\pi^{(1)}_t,
\pi^{(2)}_t\bigr)-\phi \bigl(\pi^{(1)}_s,
\pi^{(2)}_s\bigr)
\\
&&\hspace*{40pt}{}-\int_s^t\mathcal{L}\phi\bigl(
\pi^{(1)}_u,\pi^{(2)}_u\bigr)\,du-\int
_s^t \frac{\partial\phi}{\partial x}\bigl(1,
\pi_u^{(2)}\bigr)\,d\pi_u^{(3)}\biggr)
\biggr)=0,
\end{eqnarray*}
where $\psi\dvtx \tilde\Omega\to\IR$ is an arbitrary bounded,
continuous, $\mathcal{F}_s$ measurable map. Now fix such $s,t$ and
$\psi$. Then by weak convergence of $\tilde P^{(n_k)}$ to $\tilde P$
and using the moment bound in Lemma~\ref{lemEVofsquaredsuprema},
\begin{eqnarray*}
&&\lim_{k\to\infty} \tilde E^{(n_k)} \biggl( \psi(\cdot) \biggl(\phi
\bigl(\pi^{(1)}_t,\pi^{(2)}_t\bigr)-\phi
\bigl(\pi^{(1)}_s,\pi^{(2)}_s\bigr)
\\
&&\hspace*{79pt}{} -\int_s^t\mathcal{L}\phi\bigl(
\pi^{(1)}_u,\pi^{(2)}_u\bigr)\,du-\int
_s^t \frac
{\partial\phi}{\partial x}\bigl(1,
\pi_u^{(2)}\bigr)\,d\pi_u^{(3)}
\biggr) \biggr)
\\
&&\qquad =\tilde E \biggl( \psi(\cdot) \biggl(\phi\bigl(\pi^{(1)}_t,
\pi^{(2)}_t\bigr)-\phi\bigl(\pi^{(1)}_s,
\pi^{(2)}_s\bigr)
\\
&&\qquad\quad\hspace*{40pt}{} -\int_s^t\mathcal{L}\phi\bigl(
\pi^{(1)}_u,\pi^{(2)}_u\bigr)\,du-\int
_s^t \frac
{\partial\phi}{\partial x}\bigl(1,
\pi_u^{(2)}\bigr)\,d\pi_u^{(3)}
\biggr) \biggr).
\end{eqnarray*}
To complete the proof, it suffices to show that the limit on the
left-hand side above is 0. In view of the martingale property in (\ref{eq6.1}), to show this, it suffices to
prove that
for $\phi\in C^\infty_c([1,\infty)\times\IR_+)$,
%
\begin{equation}
\label{eqH4forsingletypecont} \lim_{n\to\infty} E\biggl\llvert \int
_0^t \bigl(\mathcal{L}^{(n)}\phi
\bigl(\hat X_s^{(n)},\hat Y_s^{(n)}
\bigr)-\mathcal{L}\phi\bigl(\hat X_s^{(n)},\hat
Y_s^{(n)}\bigr) \bigr)\,ds\biggr\rrvert =0
\end{equation}
and
%
\begin{equation}
\label{eqEintD^n-phi} \lim_{n\to\infty} E\biggl\llvert \int
_0^t \biggl(D^{(n)}\phi\bigl(\hat
Y_s^{(n)}\bigr)-\frac
{\partial\phi}{\partial x}\bigl(1,\hat
Y_s^{(n)}\bigr) \biggr)\,d\hat\eta^{(n)}_s
\biggr\rrvert =0.
\end{equation}
The latter is immediate upon using the smoothness of $\phi$ and the
moment estimate for $\hat\eta^{(n)}$ in (\ref{eqEVofsquaredsuprema}).
For (\ref{eqH4forsingletypecont}), we rewrite $\mathcal
{L}^{(n)}\phi$ using a Taylor expansion as follows:
\begin{eqnarray*}
&&\mathcal{L}^{(n)}\phi(x,y)
\\
&&\qquad=\lambda_1^{(n)} n^2 x \sum
_{k=0}^\infty \biggl[\frac{k-1}{n}
\,\frac
{\partial}{\partial x}\phi(x,y)+\frac{1}{2} \biggl(\frac
{k-1}{n}
\biggr)^2\,\frac{\partial^2}{\partial x^2}\phi(x,y) \biggr]\mu_1^{(n)}(k)
\\
&&\qquad\quad{} +\lambda_2^{(n)} n^2 xy \sum
_{k=0}^\infty \biggl[\frac
{k-1}{n}
\,\frac{\partial}{\partial y}\phi(x,y)+\frac{1}{2} \biggl(\frac{k-1}{n}
\biggr)^2\,\frac{\partial^2}{\partial y^2}\phi(x,y) \biggr]\mu_2^{(n)}(k)
\\
&&\qquad\quad{} +R^{(n)}(x,y)
\\
&&\qquad=c_1^{(n)}\lambda_1^{(n)} x
\,\frac{\partial}{\partial x}\phi (x,y)+\frac{1}{2}\alpha_1^{(n)}
\lambda_1^{(n)} x \,\frac{\partial^2}{\partial x^2}\phi(x,y)
\\
&&\qquad\quad{} +c_2^{(n)} \lambda_2^{(n)} x y
\,\frac{\partial}{\partial
y}\phi(x,y)+\frac{1}{2}\alpha_2^{(n)}
\lambda_2^{(n)} x y \,\frac
{\partial^2}{\partial y^2}\phi(x,y)+R^{(n)}(x,y),
\end{eqnarray*}
where the term $R^{(n)}(x,y)$ is a remainder term, which, using part
(iii) of Condition~\ref{conditioncr}, is seen to satisfy $\sup_{|x|,|y|\le L}|R^{(n)}(x,y)|\to0$ as $n\to\infty$, for any $L\in
(0,\infty)$. Furthermore, using the compact support property of $\phi
$, it follows that $\lim_{n\to\infty} E\int_0^t| R^{(n)}(\hat
X^{(n)}_s,\hat
Y^{(n)}_s)|\,ds=0$. Next note that
\begin{eqnarray*}
&&
\mathcal{L}^{(n)}\phi(x,y)-R^{(n)}(x,y)-\mathcal{L}\phi(x,y)
\\
&&\qquad=\bigl(\lambda_1^{(n)}c_1^{(n)}-
\lambda_1 c_1\bigr)x\,\frac{\partial}{\partial
x}\phi(x,y)+
\frac{1}{2}\bigl(\lambda_1^{(n)}
\alpha_1^{(n)}-\lambda_1
\alpha_1\bigr)x\,\frac{\partial^2}{\partial x^2}\phi(x,y)
\\
&&\qquad\quad{} + \bigl(\lambda_2^{(n)} c_2^{(n)}-
\lambda_2 c_2\bigr) x y \,\frac
{\partial}{\partial y}\phi(x,y)+
\frac{1}{2}\bigl(\lambda_2^{(n)}\alpha_2^{(n)}-
\lambda_2 \alpha_2\bigr)x y\,\frac{\partial^2}{\partial y^2}\phi(x,y)
\end{eqnarray*}
and therefore, in view of Condition~\ref{conditioncr},
\[
\sup_{|x|,|y|\le L}\bigl\llvert \mathcal{L}^{(n)}\phi
(x,y)-R^{(n)}(x,y)-\mathcal{L}\phi(x,y)\bigr\rrvert \to0 \qquad\mbox{as } n\to
\infty.
\]
Once more using the compact support property of $\phi$, it follows that
\[
\lim_{n\to\infty} E\int_0^t\bigl\llvert
\mathcal{L}^{(n)}\phi\bigl(\hat X_s^{(n)},\hat
Y_s^{(n)}\bigr)-R^{(n)} \bigl(\hat
X_s^{(n)},\hat Y_s^{(n)}\bigr)-
\mathcal{L}\phi \bigl(\hat X_s^{(n)},\hat
Y_s^{(n)}\bigr)\bigr\rrvert \,ds=0.
\]
Combining the above estimates, we have (\ref{eqH4forsingletypecont}),
and the result follows.
\end{pf*}

\section{\texorpdfstring{Proofs of results from Section \protect\ref{chAsymptoticBehavioroftheCatalystPopulation}}
{Proofs of results from Section 3}}
\label{secproofofAsymptoticBehavioroftheCatalystPopulation}

\subsection{\texorpdfstring{Proof of Proposition \protect\ref{propstationarydistrofX}}
{Proof of Proposition 3.1}}
\label{secproofstationarydistrofX}

Uniqueness of the invariant measure of $X$ is an immediate consequence
of the nondegeneracy of the diffusion coefficient (note that $\alpha_2 \lambda_2x\ge\alpha_2\lambda_2>0$). For existence, we will apply
an extension of the well-known Echeverria criterion for invariant
measures of Markov processes~\cite{Kurtz1991} and Theorem~5.7 of
\cite{Budhiraja2003}. This criterion, in the current context, says
that in order to establish that a probability measure $\bar\nu_1$ is
an invariant measure for $X$, it suffices to verify that for some $C\ge
0$ and all $\phi\in C^\infty_c([1,\infty))$
%
\begin{equation}
\label{eqEWK-criterion} \int_{[1,\infty)}\mathcal{L}_1
\phi(x)\bar\nu_1(dx)+C\alpha_1\lambda_1
\phi'(1)=0,
\end{equation}
where
%
\begin{equation}
\label{eq7.2} \mathcal{L}_1 \phi(x)=c_1
\lambda_1 x \phi'(x)+\tfrac{1}{2}
\alpha_1\lambda_1 x \phi''(x).
\end{equation}

We now show that (\ref{eqEWK-criterion}) holds with $\bar\nu_1=\nu_1$ and $C=\frac{p(1)}{2}$.
For $\phi\in\break C^\infty_c([1,\infty))$ and $p$ as in (\ref{eqdensityofnu1}),
\begin{eqnarray*}
&&
\int_1^\infty \biggl(c_1
\lambda_1 x \phi'(x)+\frac{1}{2}
\alpha_1 \lambda_1 x \phi''(x)
\biggr)p(x)\,dx
\\
&&\qquad=c_1\lambda_1 \theta e^{2{c_1}x/{\alpha_1}}\phi(x)
\bigg|_1^\infty-\int_1^\infty2
c_1\lambda_1 \theta\frac{c_1}{\alpha_1} e^{2{c_1}x/{\alpha_1}}
\phi(x)\,dx
\\
&&\qquad\quad{} +\frac{1}{2}\alpha_1\lambda_1 \theta
e^{2{c_1} x/{\alpha
_1}}\phi'(x) \bigg|_1^\infty-\int
_1^\infty\alpha_1
\lambda_1 \theta\frac{c_1}{\alpha_1} e^{2{c_1}x/{\alpha_1}}
\phi'(x)\,dx
\\
&&\qquad=-\frac{1}{2}\alpha_1\lambda_1 \theta
e^{2{c_1}/{\alpha
_1}}\phi'(1)=-\frac{p(1)}{2}\alpha_1
\lambda_1\phi'(1).
\end{eqnarray*}
Thus (\ref{eqEWK-criterion}) follows.

\subsection{\texorpdfstring{Proof of Theorem \protect\ref{thmtightnessandconvergenceofnu^n1}}
{Proof of Theorem 3.1}}
\label{secproofoftightnessandconvergenceofnu^n1}

Throughout this section we assume that Conditions~\ref{conditioncr},
\ref{conditionc1} and~\ref{conditionmoment} hold. This will not be
explicitly noted in the statements of the results.\vadjust{\goodbreak}

Existence of a stationary distribution $\nu^{(n)}_1$ of the
$\IS^{(n)}_X=\{\frac{l}{n}|l\in\{n,\break n+1,\ldots\}\}$ valued Markov
process $\hat X^{(n)}$ follows from the tightness of $\{\hat
X^{(n)}_t\}_{t\ge 0}$, which is a consequence of Lemma
\ref{lemsuptsupnexpmoment}. The uniqueness of the stationary
distribution follows from the irreducibility of $\hat X^{(n)}$.

In order to establish the tightness of the sequence $\{\nu^{(n)}_1\}_{n\in\IN}$, we will use the following uniform in $n$ moment
stability estimate for $\hat X^{(n)}$.

\begin{theorem}\label{thmmomentstability}
There is a $t_0\in\IR_+$ such that for all $t\geq t_0$ and $p>0$,
%
\begin{equation}
\label{eqmomentstability} \lim_{x\to\infty}\sup_{n \in\IN}
\frac{1}{x^p}E_x \bigl( \bigl(\hat X^{(n)}_{tx}
\bigr)^p \bigr)=0.
\end{equation}
\end{theorem}
\begin{pf}
Fix an $L>1$, and let $\tau^{(n)}:=\inf\{t\dvtx  \hat X^{(n)}_t\le L\}$.
Observe that if $t\in[(N-1)\rho,N\rho)$ for some $N\in\IN$, then,
following arguments as in the proof of Lemma \ref
{lemsuptsupnexpmoment}, for $x>L$,
\[
P_x\bigl(\tau^{(n)}>t\bigr)\le P_x \Biggl(
\sum_{j=1}^{N-1}\bigl(v^{(n)}_j-
\delta \kappa\rho\bigr)>\delta(L-x) \Biggr)\le e^{\delta(x-L-\delta\kappa\rho(N-1))}.
\]
Thus we have that
\[
\sup_{n\in\IN}P_x\bigl(\tau^{(n)}>t\bigr)\le
\gamma_1 e^{\delta x} e^{-\gamma_2 t},
\]
where $\gamma_i\in(0,\infty)$, $i=1,2$. The above estimate along
with Lemma~\ref{lemsuptsupnexpmoment} implies, for $n\in\IN$,
\begin{eqnarray*}
E_x e^{{\delta}\hat X^{(n)}_t/{2}}&=& E_x \bigl( 1_{\{\tau
^{(n)}\le t\}}
e^{{\delta}\hat X^{(n)}_t/{2}} \bigr)+E_x \bigl( 1_{\{\tau^{(n)}> t\}}
e^{{\delta}\hat X^{(n)}_t/{2}} \bigr)
\\
&\le&\tilde d e^{ \delta L}+ \bigl(\gamma_1 e^{\delta x}e^{-\gamma_2
t}
\bigr)^{1/2} \bigl( E_x \bigl(e^{\delta\hat
X^{(n)}_t}
\bigr) \bigr)^{1/2}
\\
&\le&\tilde de^{ \delta L}+ \bigl(\gamma_1 e^{\delta x}e^{-\gamma_2 t}
\bigr)^{1/2}\bigl(\tilde de^{\delta x}\bigr)^{1/2} \le
d_1\bigl(1+e^{\delta x} e^{-{\gamma_2}/{2} t}\bigr),
\end{eqnarray*}
where $\tilde d$ is as in Lemma~\ref{lemsuptsupnexpmoment} and
$d_1\in(0,\infty)$ is some constant, independent of $n$. Fix $p>0$.
Then, for some $d_2\in(0,\infty)$, we have
\[
\sup_{n \in\IN}\frac{E_x(\hat X^{(n)}_{tx})^p}{x^p}\le\sup_{n \in
\IN}\frac{d_2 E_x e^{{\delta}\hat X^{(n)}_{tx}/{2}}}{x^p}
\le \sup_{n \in\IN} \frac{d_1 d_2(1+e^{\delta x} e^{-{\gamma_2
}tx/{2}}) }{x^p}.
\]
Choose $t_0$ large enough such that $\frac{\gamma_2}{2} t_0 >\delta
$. Then for $t\ge t_0$
\[
\lim_{x\to\infty}\sup_{n\in\IN}\frac{E_x ((\hat
X^{(n)}_{tx})^p )}{x^p}=0.
\]
The result follows.
\end{pf}

As a consequence of Theorem~\ref{thmmomentstability}, we have the
following result.
For $\delta\in(0,\infty)$, define the return time to a compact set
$C\subset[1,\infty)$ by $\tau_C^{(n)}(\delta):=\inf\{t\geq\delta
|\hat X^{(n)}_{t}\in C\}$.\vadjust{\goodbreak}
%
\begin{theorem}\label{thmintegratingtillreturntime}
There are $\tilde c,\hat\delta\in(0,\infty)$ and a compact set
$C\subset[1,\infty)$ such that
\[
\sup_n E_x \biggl(\int_0^{\tau_C^{(n)}(\hat\delta)}
\ln\bigl(\hat X^{(n)}_{t}\bigr)\,dt \biggr)\leq\tilde c
x^3,\qquad x\geq1.
\]
\end{theorem}
\begin{pf}
Note that $\ln(\hat X^{(n)}_t)\ge0$ since $\hat X^{(n)}_t\ge1$.
Applying Theorem~\ref{thmmomentstability} with $p=3$, we have that
there is an $L\in(1,\infty)$ such that with $C:=\{x\in\IR_+ | x\leq
L\}$, for all $x\in C^c$,
%
\begin{equation}
\label{eqexpmomentoutsidecompact} \sup_n E_x \bigl(
\bigl(\hat X^{(n)}_{t_0 x} \bigr)^3 \bigr)\leq
\frac{1}{2} x^3,
\end{equation}
where $t_0$ is as in Theorem~\ref{thmmomentstability}. Let $\hat
\delta:=t_0 L$ and
\[
\tau^{(n)}:=\tau_C^{(n)}(\hat\delta) =\inf\bigl
\{t\geq\hat\delta |\hat X^{(n)}_{t}\leq L\bigr\}.
\]
Consider a sequence of stopping times defined as follows:
\[
\sigma^{(n)}_0:=0,\qquad \sigma^{(n)}_m:=
\sigma^{(n)}_{m-1}+t_0 \bigl(\hat
X^{(n)}_{\sigma^{(n)}_{m-1}}\vee L \bigr),\qquad m\in\IN.
\]
Let $m_0^{(n)}=\min\{m\geq1 |\hat X^{(n)}_{\sigma^{(n)}_m}\leq L\}$, and
%
\begin{equation}
\label{lyapfn} \hat V^{(n)}(x):=E_x \biggl(\int
_0^{\tau^{(n)}}\ln\bigl(\hat X^{(n)}_t
\bigr)\,dt \biggr).
\end{equation}
Then
\[
\hat V^{(n)}(x)\leq E_x \biggl(\int_0^{\sigma^{(n)}_{m_0^{(n)}}}
\ln \bigl(\hat X^{(n)}_t\bigr)\,dt \biggr)=\sum
_{k=0}^\infty E_x \biggl(\int
_{\sigma
^{(n)}_k}^{\sigma^{(n)}_{k+1}}\ln\bigl(\hat X^{(n)}_t
\bigr)\,dt 1_{\{k<m_0^{(n)}\}
} \biggr).
\]
Let $\mathcal{F}^{(n)}_t:=\sigma\{\hat X^{(n)}_s |0\leq s\leq t\}$.
We claim that there is a $c_0\in(0,\infty)$ such that for all $n,k\in
\IN$, $x\geq1$
%
\begin{equation}
\label{eqclaimaa} E_x \biggl(\int_{\sigma^{(n)}_k}^{\sigma^{(n)}_{k+1}}
\ln\bigl(\hat X^{(n)}_t\bigr)\,dt\Big|\mathcal{F}^{(n)}_{\sigma^{(n)}_k}
\biggr) 1_{\{
k<m_0^{(n)}\}} \leq c_0 \bigl(\hat X^{(n)}_{\sigma^{(n)}_k}
\bigr)^3 1_{\{k<m_0^{(n)}\}}.
\end{equation}
Due to the strong Markov property, to prove the claim
it suffices to show that for some $c_0\in(0,\infty)$ and for all
$n\in\IN, x\ge1$
\[
E_x \biggl(\int_0^{\sigma^{(n)}_1}\ln\bigl(
\hat X^{(n)}_t\bigr)\,dt \biggr)\le c_0
x^{3}.
\]
Note that for $x\ge1$, $\sigma^{(n)}_1= t_0(x\vee L)\le\tilde c_0
x$, where $\tilde c_0=t_0L$. Using this bound along with Lemma \ref
{lemEVofsquaredsuprema}, we get, for some $\hat c_0\in(0,\infty)$,
\[
E_x \Bigl(\sup_{ t\le\sigma^{(n)}_1}\ln\bigl(\hat X^{(n)}_t
\bigr) \Bigr)\le \ln \Bigl(E_x \Bigl(\sup_{ t\le\sigma^{(n)}_1}\hat
X^{(n)}_t \Bigr) \Bigr)\le\ln \bigl(x^2
e^{K (\tilde c_0 x)^2 } \bigr)\le\hat c_0 x^2.
\]
The claim follows.\vadjust{\goodbreak}

From the estimate (\ref{eqclaimaa}), we now have
%
\begin{equation}
\label{eq7.7} \sup_n \hat V^{(n)}(x)\leq c_0
\sup_n E_x \Biggl(\sum_{k=0}^{m_0^{(n)}-1}
\bigl(\hat X^{(n)}_{\sigma^{(n)}_k} \bigr)^3 \Biggr).
\end{equation}
Note that $ \{\hat X^{(n)}_{\sigma^{(n)}_k} \}_{k\in\IN_0}$
is a Markov chain with transition probability kernel
\[
\check P^{(n)}(x,A):= P^{(n)}_{t_0(x\vee L)}(x,A),\qquad x\in[1,
\infty ), A\in\mathcal{B}\bigl([1,\infty)\bigr),
\]
where $P^{(n)}_t$ is the transition probability kernel for $\hat
X^{(n)}$. Using (\ref{eqexpmomentoutsidecompact}) and Lem\-ma~\ref
{lemsuptsupnexpmoment}, we get that for any $L\in(1,\infty)$
there exists a $\tilde b\in(0,\infty)$ such that for all $x\in
[1,\infty)$
%
\begin{eqnarray}
\label{eqdriftcondition}
\sup_n \int_1^\infty y^3
\check P^{(n)}(x,dy)&=&\sup_n \int_1^\infty
y^3 P^{(n)}_{t_0(x\vee L)}(x,dy)
\nonumber
\\
&=&\sup_n \int_1^\infty
y^3 P^{(n)}_{t_0x}(x,dy) 1_{\{x>L\}}\nonumber\\[-8pt]\\[-8pt]
&&{}+
\sup_n \int_1^\infty y^3
P^{(n)}_{t_0 L}(x,dy) 1_{\{x\leq L\}}
\nonumber\\
&\leq& x^3-\frac{1}{2}x^3+\tilde b
1_{[0,L]}(x).
\nonumber
\end{eqnarray}
The above inequality along with Theorem 14.2.2 of~\cite{Meyn1993} yields
\begin{eqnarray*}
\sup_n E_x \Biggl(\sum_{k=0}^{m_0^{(n)}-1}
\bigl(\hat X^{(n)}_{\sigma
^{(n)}_k} \bigr)^3 \Biggr)&\leq& 2
\Biggl(x^3+\sup_n E_x \Biggl(\sum
_{k=0}^{m_0^{(n)}-1} \tilde b 1_{[0,L]} \bigl(\hat
X^{(n)}_{\sigma
^{(n)}_k} \bigr) \Biggr) \Biggr)
\\
&=& 2 \bigl(x^3+\tilde b 1_{[0,L]}(x) \bigr)\le\tilde c
x^3,
\end{eqnarray*}
where the equality in the last display follows from the fact that $\hat
X^{(n)}_{\sigma^{(n)}_k}>L$ for $1\leq k<m_0^{(n)}$.
The result follows now on combining the last estimate with~(\ref
{eq7.7}).
\end{pf}

The following theorem is proved exactly as Proposition 5.4 of \cite
{Dai1995}. The proof is omitted.
%
\begin{theorem}\label{thmestimatefX^nbyV^n}
Let $f\dvtx [1,\infty)\to\IR_+$ be a measurable function. Define for
$\hat\delta\in(0,\infty)$ and a compact set $C\subset[1,\infty)$
\[
V^{(n)}(x):=E_x \biggl(\int_0^{\tau^{(n)}_C(\hat\delta)}f
\bigl(\hat X^{(n)}_t\bigr)\,dt \biggr),\qquad x\in[1,\infty).
\]
If $\sup_{n\in\IN} V^{(n)}$ is everywhere finite and uniformly
bounded on $C$, then there exists a $\hat\kappa\in(0,\infty)$ such
that for all $n\in\IN, t>0, x\in[1,\infty)$
\[
\frac{1}{t}E_x \bigl(V^{(n)}\bigl(\hat
X^{(n)}_t\bigr) \bigr)+\frac{1}{t}\int
_0^t E_x \bigl(f\bigl(\hat
X^{(n)}_s\bigr) \bigr)\,ds\le\frac
{1}{t}V^{(n)}(x)+
\hat\kappa.
\]
\end{theorem}

We now return to the proof of Theorem~\ref{thmtightnessandconvergenceofnu^n1}
and establish the tightness of $\{\nu^{(n)}_1\}_{n\in\IN}$.
We will apply
Theorem~\ref{thmestimatefX^nbyV^n} with $f(x):=\ln(x)$, and
$\hat\delta, C$ as in Theorem~\ref{thmintegratingtillreturntime}.
Since $\nu^{(n)}_1$ is an invariant measure for $\hat X^{(n)}$, we
have for nonnegative, real valued, measurable functions $\Phi$ on
$[1,\infty)$
%
\begin{equation}
\label{eq7.9} \int_1^\infty E_x
\bigl(\Phi\bigl(\hat X^{(n)}_t\bigr) \bigr)
\nu^{(n)}_1(dx)=\int_1^\infty
\Phi(x)\nu^{(n)}_1(dx).
\end{equation}
Fix $k\in\IN$ and let $V^{(n)}_k(x):=V^{(n)}(x)\wedge k$. Let
\[
\Psi^{(n)}_k(x):=\frac{1}{t}V^{(n)}_k(x)-
\frac{1}{t}E_x \bigl(V^{(n)}_k\bigl(\hat
X^{(n)}_t\bigr) \bigr).
\]
By (\ref{eq7.9}), we
have that $\int_1^\infty\Psi^{(n)}_k(x)\nu^{(n)}_1(dx)=0$. Let
\[
\Psi^{(n)}(x):=\frac{1}{t}V^{(n)}(x)-
\frac{1}{t}E_x \bigl(V^{(n)}\bigl(\hat
X^{(n)}_t\bigr) \bigr).
\]
By the monotone convergence theorem $\Psi^{(n)}_k(x)\to\Psi^{(n)}(x)$ as $k\to\infty$.
We next show that $\Psi^{(n)}_k(x)$ is bounded from below for all
$x\in[1,\infty)$: if $V^{(n)}(x)\le k$, then
\begin{eqnarray*}
\Psi^{(n)}_k(x)&=&\frac{1}{t}V^{(n)}_k(x)-
\frac{1}{t} E_x \bigl(V^{(n)}_k\bigl(
\hat X^{(n)}_t\bigr) \bigr)
\\
&\ge&\frac{1}{t} V^{(n)}(x)-\frac{1}{t} E_x
\bigl(V^{(n)}\bigl(\hat X^{(n)}_t\bigr) \bigr)\ge-
\hat\kappa,
\end{eqnarray*}
where the last inequality follows from Theorem \ref
{thmestimatefX^nbyV^n}. If $V^{(n)}(x)\ge k$,
\[
\Psi^{(n)}_k(x)=\frac{1}{t}k-\frac{1}{t}
E_x \bigl(V^{(n)}_k\bigl(\hat
X^{(n)}_t\bigr) \bigr)\ge0.
\]
Thus $\Psi^{(n)}_k(x)\ge-\hat\kappa$ for all $x\ge1$. By Fatou's
lemma, we have
\[
\int_1^\infty\Psi^{(n)}(x)
\nu^{(n)}_1(dx)\le\liminf_{k\to\infty
}\int
_1^\infty\Psi^{(n)}_k(x)
\nu^{(n)}_1(dx)=0.
\]
By Theorem~\ref{thmestimatefX^nbyV^n} we have $\Psi^{(n)}(x)\ge
\frac{1}{t}\int_0^t E_x (f(\hat X^{(n)}_s) )\,ds-\hat\kappa
$. Combining this with the last display, we have
\[
0\ge\int_1^\infty\Psi^{(n)}(x)
\nu^{(n)}_1(dx)\ge\frac{1}{t}\int
_0^t\int_1^\infty
E_x \bigl(f\bigl(\hat X^{(n)}_s\bigr) \bigr)
\nu^{(n)}_1(dx)\,ds-\hat\kappa.
\]
Using the invariance property of $\nu^{(n)}_1$ once more, we see that
the first term on the right-hand side above equals $\int f(x) \nu^{(n)}_1(dx)$, and
therefore\break
$\int f(x) \nu^{(n)}_1(dx) \le\hat\kappa$.
This completes the proof of tightness.

The tightness of $\{\nu^{(n)}_1\}_{n\in\IN}$ implies that every
subsequence of $\{\nu^{(n)}_1\}$ has a convergent subsequence. Call
such a limit $\nu^*_1$. Theorem~\ref{thmcrdiffusion} and the
stationarity of $\nu^{(n)}_1$ imply that $\nu^*_1$ is a stationary
distribution of $X$. Since the stationary distribution of $X$ is
unique, we have $\nu^*_1=\nu_1$, which completes the
proof.\looseness=-1

\section{\texorpdfstring{Proofs of Theorems \protect\ref{thmfastcatalystweaklimit} and
\protect\ref{fastcatalystBPaveraging}}
{Proofs of Theorems 4.1 and 4.2}}
\label{secDiffusionLimitoftheReactantunderFastCatalystDynamics}

\subsection{\texorpdfstring{Proof of Theorem \protect\ref{thmfastcatalystweaklimit}}
{Proof of Theorem 4.1}}

In order to prove the result, we will verify that the assumptions of
Theorem II.1 (more precisely, those in the remark following Theorem
II.1) in~\cite{Skorokhod1989}, pages 78\vspace*{1pt} and 79, hold. For this, it
suffices to show that for all $k\in\IN$, $\Phi\in\operatorname
{BM}(\IR_+^k)$, $\phi\in C^\infty_c(\IR_+)$ and $0\leq
t_1<t_2<\cdots<t_{k+1}<T<\infty$, there exists a sequence $h_n$ with
$\lim_{n\to\infty} h_n=0$ and
%
\begin{eqnarray}
\label{eqstochaveragcriterion}
&&
\sup_{t\in[t_{k+1},T]}\bigl\llvert E \bigl[\Phi
\bigl(\check Y^{(n)}_{t_1},\ldots,\check Y^{(n)}_{t_k}
\bigr) \bigl(\phi \bigl(\check Y^{(n)}_{t+h_n} \bigr)-\phi
\bigl(\check Y^{(n)}_{t} \bigr)-h_n\check L\phi
\bigl(\check Y^{(n)}_t \bigr) \bigr) \bigr]\bigr\rrvert
\nonumber\\[-8pt]\\[-8pt]
&&\qquad =o(h_n),\nonumber
\end{eqnarray}
where $\check{\mathcal{L}}$ is given as
%
\begin{equation}
\label{eq8.2} \check{\mathcal{L}}\phi(y):= c_2\lambda_2
m_X y\phi'(y)+\tfrac
{1}{2}\alpha_2
\lambda_2 m_X y\phi''(y),\qquad \phi
\in C^\infty_c(\IR_+).
\end{equation}

Letting $X^\circ_t:=\check X^{(n)}_{t/a_n}$, $t\ge0$, we see, using
scaling properties of the Skorohod map and straight forward martingale
characterization results, that $X^\circ$ has the same probability law
as the process $X$ that was introduced in Proposition \ref
{propuniqunessofX,Y} with initial value $X_0=x_0$. The following
uniform moment bound will be used in the proof of Theorem \ref
{thmfastcatalystweaklimit}.

\begin{lemma}\label{lemunifmombdX}
There exists a $\delta_0 \in(0, \infty)$, such that whenever $X$ is
as in Proposition~\ref{propuniqunessofX,Y} with initial value
$X_0=x$, for some $x\in[1,\infty)$, we have
\[
\sup_{0\le t<\infty}E_x \bigl( e^{\delta_0 X_t} \bigr)=:d(
\delta_0,x)<\infty.
\]
\end{lemma}
\begin{pf}
We begin by establishing exponential moment estimates for the increase
of $ X$ over time intervals of length $l\rho$ when the process is away
from the boundary 1, where $\rho>0$ and $l\geq1$.
Fix $\rho\in(0,\infty)$. Let $a:=c_1\lambda_1$,
$b:=\alpha_1\lambda_1$ and $\delta\in(0,-\frac{a}{b}\wedge1)$.
Note that\vspace*{1pt} in view of Condition~\ref{conditionc1}, $a<0$.
Define $\sigma_r:=\inf\{t\in[0,\infty)|\int_0^t X_s \,ds>r\}$ and
$\rho_{l,r}:=l\rho\wedge\sigma_r$.
Then
\begin{eqnarray*}
&&
E_x \biggl(\exp \biggl(\delta a\int_0^{\rho_{l,r}}
X_s \,ds + \delta \sqrt{b}\int_0^{\rho_{l,r}}
\sqrt{ X_s}\,dB^X_s \biggr) \biggr)
\\
&&\qquad=E_x \biggl(\exp \biggl(\bigl(\delta a+\delta^2 b
\bigr)\int_0^{\rho_{l,r}} X_s \,ds\\
&&\qquad\hspace*{52pt}{}+ \delta
\sqrt{b}\int_0^{\rho_{l,r}} \sqrt{ X_s}\,dB^X_s-
\delta^2 b\int_0^{\rho_{l,r}}
X_s \,ds \biggr) \biggr)
\\
&&\qquad\le\sqrt{ E_x
e^{2\rho_{l,r}(\delta a+\delta^2 b)}} \\
&&\qquad\quad{}\times\sqrt{E_x\exp \biggl(2\delta
\sqrt{b}\int_0^{\rho_{l,r}} \sqrt{ X_s}\,dB^X_s-
\frac{(2\delta\sqrt{b})^2}{2}\int_0^{\rho_{l,r}} X_s
\,ds \biggr)},
\end{eqnarray*}
where the inequality follows on noting that $\rho_{l,r}\le\int_0^{\rho_{l,r}} X_s \,ds$ and $\delta\in(0,-\frac{a}{b})$. Using the
super martingale property of the stochastic exponential, we have that
the second term on the right-hand side of the last display is bounded
by 1.
Thus, sending $r\to\infty$, we have, with $-\theta:=\delta a+\delta^2 b<0$,
%
\begin{equation}
\label{eq8.3}\qquad E_x \biggl(\exp \biggl(\delta a\int
_0^{l\rho} X_s \,ds + \delta\sqrt {b}
\int_0^{l\rho} \sqrt{ X_s}\,dB^X_s
\biggr) \biggr)\le e^{l\rho
(\delta a+\delta^2 b)}=e^{-\theta l\rho}.
\end{equation}

Next, for $x\in[1,\infty)$, we have by application of It\^{o}'s
formula that for $t\le\rho$ and $\tilde\delta\le\delta$,
%
\begin{equation}
\label{eq8.4} E_x \bigl(e^{\tilde\delta X_t} \bigr)\le
e^{\tilde\delta x}+\tilde \delta e^{\tilde\delta}E_x
\eta_\rho\le e^{\tilde\delta x}+x C_1(\rho,\delta),
\end{equation}
where $C_1(\rho,\delta)\in(0,\infty)$ and the last inequality
follows by an application of Gronwall's lemma and the Lipschitz
property of the Skorohod map; see~(\ref{eqSkorohodLipschitz}).

Using the above estimates, we will now establish certain uniform
estimates on the tail probabilities of $X_{k\rho}$, which will lead to
exponential moment estimates at these time points. Fix $L>1$, and let
\[
\tau_j:=\inf\bigl\{t\geq(j-1)\rho|X_t\leq L\bigr\}
\wedge j\rho\quad\mbox{and}\quad e_j:=X_{j\rho}-X_{\tau_j},\qquad
j\geq1,
\]
$e_0=0$. Fix $k\in\IN$,
and let
\[
M:= \max \Bigl\{j=1,\ldots,k \big|\inf_{(j-1)\rho\leq s\leq
j\rho} X_s\leq L \Bigr\},
\]
if there is an $s\in[0,k\rho]$ such that $X_s\leq L$, and set $M$
equal to 0 otherwise.
Let
\[
v_j:=\int_{(j-1)\rho}^{j\rho}a
X_s \,ds+ \int_{(j-1)\rho}^{j\rho
}
\sqrt{bX_s}\,dB^X_s,\qquad j\geq1.
\]
Then $ X_{k\rho} = X_{M\rho}+\sum_{j=M+1}^k v_j$.
Letting $\zeta_i:=e_i+\sum_{j=i+1}^k v_j$, we have, using~(\ref
{eq8.3}),
\begin{eqnarray*}
P_x( X_{k\rho}>K)&\leq& P_x
\Biggl(X_{M\rho}+\sum_{j=M+1}^k
v_j>K \Biggr)\\
&\leq& P_x \Bigl(\max_{0\leq i\leq k}
\zeta_i>K-L \Bigr)
\\
&\leq&\sum_{i=0}^k P_x(
\zeta_i>K-L) \leq\sum_{i=0}^k
E_x \bigl(e^{{\delta}\zeta_i/{2}} \bigr)e^{-{\delta}(K-L)/{2}}
\\
&\leq&\sum_{i=0}^k \bigl(E_x
e^{\delta e_i} \bigr)^{1/2} \bigl(E_x
e^{\delta\sum_{j=i+1}^k v_j} \bigr)^{1/2}e^{-
{\delta}(K-L)/{2}}
\\
&\leq&\sum_{i=0}^k \bigl(E_x
e^{\delta e_i} \bigr)^{1/2} e^{-(k-i)\theta\rho/{2}}
e^{-{\delta}(K-L)/{2}}.
\end{eqnarray*}
Next note that, from (\ref{eq8.4}),
\begin{eqnarray*}
E_x e^{\delta e_i}&\le& E_x \bigl(e^{\delta[X_{i\rho}-X_{\tau_i}]}
1_{\tau_i<i\rho} \bigr)+1=E_x \bigl( e^{-\delta X_{\tau_i}}1_{\tau
_i<i\rho}
E_{X_{\tau_i}} \bigl(e^{\delta X_{i\rho}} \bigr) \bigr)+1
\\
&\le& E_x \bigl[e^{-\delta X_{\tau_i}}1_{\tau_i<i\rho}
\bigl(e^{\delta X_{\tau_i}}+X_{\tau_i} C_1(\rho,\delta) \bigr)
\bigr]+1\equiv C_2(\rho,\delta)+1.
\end{eqnarray*}
Hence,
\begin{eqnarray*}
P_x( X_{k\rho}>K)&\leq&\bigl( C_2(\rho,
\delta)+1\bigr)^{1/2} e^{-
{\delta}(K-L)/{2}}\sum
_{l=0}^k e^{-l\rho\theta/{2}}
\\
&\leq&\bigl( C_2(\rho,\delta)+1\bigr)^{1/2}
\frac{ e^{-{\delta
}(K-L)/{2}}}{1-e^{-\rho\theta/{2}}}.
\end{eqnarray*}
The last estimate yields, analogously to (\ref{eq5.38}),
%
\begin{equation}
\label{eq8.5} \sup_{k\in\IN_0}E_x\bigl(e^{{\delta} X_{k\rho}/4}\bigr)
\leq Ce^{
{\delta}x/{2}}\qquad\mbox{for all }x\in[1,\infty)
\end{equation}
for some $C\in(0,\infty)$.
Finally, letting $\delta_0:=\frac{\delta}{4}$, we have from (\ref
{eq8.4})
for $t\in((k-1)\rho,k\rho]$, $k\ge1$,
\begin{eqnarray*}
E_x \bigl(e^{\delta_0 X_t}\bigr)&=&E_x
\bigl(E_{X_{(k-1)\rho}} \bigl(e^{\delta_0
X_t}\bigr) \bigr)\\
&\le& E_x
\bigl(e^{\delta_0 X_{(k-1)\rho}}+ X_{(k-1)\rho} C_1(\rho,\delta) \bigr)
\\
&\le& C e^{\frac{\delta}{2}x} \biggl(1+\frac{1}{\delta_0}C_1(\rho,
\delta ) \biggr).
\end{eqnarray*}
The result follows.
\end{pf}

\begin{remark}\label{remunifmombdcheckX^n}
Note that Lemma~\ref{lemunifmombdX} and the scaling property noted
above that lemma say that for all $x\in[1,\infty)$
\[
\sup_{n\in\IN}\sup_{0\le t<\infty}E_x \bigl(e^{\delta_0 \check
X^{(n)}_t}
\bigr)<\infty.
\]
\end{remark}

We now prove Theorem~\ref{thmfastcatalystweaklimit} by showing
(\ref{eqstochaveragcriterion}).
Let, for $\phi\in C^\infty_c(\IR_+)$,
%
\begin{equation}
\label{eq8.6}\qquad \mathcal{L}_{x}\phi(y):= c_2
\lambda_2 x y\phi'(y)+\tfrac
{1}{2}
\alpha_2\lambda_2 x y\phi''(y),\qquad
(x,y)\in[1,\infty )\times\IR_+.
\end{equation}
Then
%
\begin{eqnarray}
\label{eqab}
&&
E \bigl[\Phi \bigl(\check Y^{(n)}_{t_1},\ldots,
\check Y^{(n)}_{t_k} \bigr) \bigl(\phi\bigl(\check
Y^{(n)}_{t+h_n}\bigr)-\phi\bigl(\check Y^{(n)}_{t}
\bigr) \bigr) \bigr]
\nonumber
\\
&&\qquad=E \biggl[\Phi \bigl(\check Y^{(n)}_{t_1},\ldots,\check
Y^{(n)}_{t_k} \bigr) \int_t^{t+h_n}
\mathcal{L}_{\check X^{(n)}_s}\phi\bigl(\check Y^{(n)}_t\bigr)\,ds
\biggr]
\\
&&\qquad\quad{} +E \biggl[\Phi \bigl(\check Y^{(n)}_{t_1},\ldots,\check
Y^{(n)}_{t_k} \bigr) \int_t^{t+h_n}
\bigl(\mathcal{L}_{\check
X^{(n)}_s}\phi\bigl(\check Y^{(n)}_s
\bigr)- \mathcal{L}_{\check X^{(n)}_s}\phi \bigl(\check Y^{(n)}_t
\bigr) \bigr)\,ds \biggr].\hspace*{-15pt}\nonumber
\end{eqnarray}
For the second term, we have, using Remark~\ref{remunifmombdcheckX^n}
and the fact that the function $\Phi$ is bounded and $\phi$ as
well as its derivatives are continuous with bounded support, that
%
\begin{eqnarray}
\label{eqsecondterm} &&
\sup_{t\in[t_{k+1},T]} E \biggl|\Phi \bigl(\check
Y^{(n)}_{t_1},\ldots,\check Y^{(n)}_{t_k}
\bigr) \int_t^{t+h_n} \bigl(\mathcal{L}_{\check X^{(n)}_s}
\phi\bigl(\check Y^{(n)}_s\bigr)- \mathcal{L}_{\check X^{(n)}_s}
\phi\bigl(\check Y^{(n)}_t\bigr) \bigr)\,ds \biggr|
\nonumber
\\
&&\qquad= \sup_{t\in[t_{k+1},T]} E \biggl|\Phi \bigl(\check Y^{(n)}_{t_1},
\ldots,\check Y^{(n)}_{t_k} \bigr) \nonumber\\
&&\qquad\quad\hspace*{48pt}{}\times\int_t^{t+h_n}
\biggl[c_2\lambda_2\check X^{(n)}_s
\bigl(\check Y^{(n)}_s\phi'\bigl( \check
Y^{(n)}_s\bigr)-\check Y^{(n)}_t
\phi'\bigl(\check Y^{(n)}_t\bigr) \bigr)
\\
&&\hspace*{93.3pt}\qquad\quad{} +\frac{1}{2}\alpha_2\lambda_2 \check
X^{(n)}_s \bigl( \check Y^{(n)}_s
\phi''\bigl(\check Y^{(n)}_s
\bigr)-\check Y^{(n)}_t\phi''
\bigl(\check Y^{(n)}_t\bigr) \bigr) \biggr]\,ds
\biggr|\nonumber\\
&&\qquad=o(h_n).\nonumber
\end{eqnarray}
Recalling the definition of $X^\circ$ above Lemma~\ref{lemunifmombdX},
the first expected value on the right-hand side in (\ref
{eqab}) equals
\[
h_nE \biggl[\Phi \bigl(\check Y^{(n)}_{t_1},\ldots,\check Y^{(n)}_{t_k} \bigr) \frac{1}{h_n a_n}\int
_{t a_n}^{t a_n+h_n a_n}\mathcal{L}_{ X^\circ
_s}\phi\bigl(
\check Y^{(n)}_t\bigr)\,ds \biggr].
\]
Thus
\begin{eqnarray*}
&&
E \bigl[\Phi \bigl(\check Y^{(n)}_{t_1},\ldots,\check
Y^{(n)}_{t_k} \bigr) \bigl(\phi\bigl(\check
Y^{(n)}_{t+h_n}\bigr)-\phi\bigl(\check Y^{(n)}_{t}
\bigr)-h_n\check{\mathcal{L}}\phi\bigl(\check Y^{(n)}_t
\bigr) \bigr) \bigr]
\\
&&\qquad=E \biggl[\Phi \bigl(\check Y^{(n)}_{t_1},\ldots,\check
Y^{(n)}_{t_k} \bigr) h_n \biggl(\frac{1}{h_n a_n}
\int_{ta_n}^{ta_n+h_na_n} \mathcal{L}_{ X^\circ_s}\phi
\bigl(\check Y^{(n)}_t\bigr)\,ds-\check{\mathcal{L}}\phi\bigl(
\check Y^{(n)}_t\bigr) \biggr) \biggr]\\
&&\qquad\quad{}+o(h_n)
\\
&&\qquad=E \biggl[\Phi \bigl(\check Y^{(n)}_{t_1},\ldots,\check
Y^{(n)}_{t_k} \bigr) h_n \biggl(c_2
\lambda_2\check Y^{(n)}_t\phi'
\bigl(\check Y^{(n)}_t\bigr)+\frac
{1}{2}
\alpha_2\lambda_2 \check Y^{(n)}_t
\phi''\bigl(\check Y^{(n)}_t
\bigr) \biggr)
\\
&&\hspace*{173.4pt}{}\times \biggl(\frac{1}{h_na_n}\int_{ta_n}^{ta_n+h_na_n}
X^\circ_s \,ds - m_X \biggr)
\biggr]\\
&&\qquad\quad{}+o(h_n).
\end{eqnarray*}
To complete the proof, it thus remains to show that for some sequence
$\{h_n\}$ with $\lim_{n\to\infty} h_n=0$
%
\begin{equation}
\label{eqc}\quad E\biggl\llvert \frac{1}{h_na_n}\int_{ta_n}^{ta_n+h_na_n}
X^\circ_{s}\,ds-m_X\biggr\rrvert =E\biggl\llvert
\frac{1}{h_na_n}\int_{ta_n}^{ta_n+h_na_n}
X_{s}\,ds-m_X\biggr\rrvert
\end{equation}
converges to 0 uniformly in $t\in[t_{k+1},T]$.
From the ergodicity of $X$ and the moment estimate in Lemma \ref
{lemunifmombdX}, it follows that
%
\begin{equation}
\label{eqb} E\biggl\llvert \frac{1}{ta_n}\int_0^{ta_n}
X_{s}\,ds-m_X\biggr\rrvert \to0\qquad \mbox{as }n\to\infty.
\end{equation}
The above result, along with Lemma~\ref{lemshiftintegrallimits},
below, implies that there is a sequence $\{h_n\}$ such that $\lim_{n\to
\infty} h_n= 0$,
and the expression in (\ref{eqc}) converges to 0 uniformly in $t\in
[t_{k+1},T]$. This completes the proof.

The proof of the following lemma is adapted from Lemma II.9, page 137,
in~\cite{Skorokhod1989}.
%
\begin{lemma}\label{lemshiftintegrallimits}
Let $0\le t_{k+1}<T<\infty$ and $a_n\to\infty$ monotonically as
$n\to\infty$. If for all $t\in[t_{k+1},T]$
\[
E\biggl\llvert \frac{1}{ta_n}\int_{0}^{ta_n}X_s\,ds-m_X
\biggr\rrvert \to0\qquad \mbox{as }n\to\infty,
\]
then there is a sequence $\{h_n\}$ such that $h_n\to0$ as $n\to\infty
$, and
\[
\sup_{t\in[t_{k+1},T]}E\biggl\llvert \frac{1}{h_na_n}\int_{ta_n}^{ta_n+h_na_n}X_s\,ds-m_X
\biggr\rrvert \to0\qquad \mbox{as }n\to \infty.
\]
\end{lemma}
\begin{pf}
Let
$\alpha(\tau):=\sup_{u>\tau}E\llvert \frac{1}{u}\int_{0}^{u}X_s\,ds-m_X\rrvert $. Note that $\alpha(\tau)$ converges
monotonically to 0 as $\tau\to\infty$.
For $t\in[t_{k+1},T]$ we have
\begin{eqnarray*}
\hspace*{-4pt}&&
E\biggl\llvert \frac{1}{h_na_n}\int_{ta_n}^{ta_n+h_na_n}X_s\,ds-m_X
\biggr\rrvert
\\
\hspace*{-4pt}&&\quad=E\biggl\llvert \frac{ta_n+h_na_n}{h_na_n}\frac{1}{ta_n+h_na_n}\int_{0}^{ta_n+h_na_n}
X_s\,ds-\frac{ta_n}{h_na_n}\frac
{1}{ta_n}\int_{0}^{ta_n}X_s\,ds-m_X
\biggr\rrvert
\\
\hspace*{-4pt}&&\quad\leq\frac{ta_n+h_na_n}{h_na_n} \alpha(ta_n+h_na_n)+
\frac
{ta_n}{h_na_n} \alpha(ta_n)\leq\frac{3T}{h_n}
\alpha(t_{k+1}a_n)
\end{eqnarray*}
for all $n$ such that $h_n\le T$. Note that the right-hand side of the
last display is independent of $t\in[t_{k+1},T]$. Choosing $h_n=\sqrt {\alpha(t_{k+1}a_n)}$, the lemma follows.
\end{pf}

\subsection{\texorpdfstring{Proof of Theorem \protect\ref{fastcatalystBPaveraging}}
{Proof of Theorem 4.2}}
\label{secproofoffastcatalystweaklimit}

As in the proof of Theorem~\ref{thmfastcatalystweaklimit}, it
suffices to show that for all $k\in\IN$, $\Phi\in\operatorname
{BM}(\IR_+^k)$, $\phi\in C^\infty_c(\IR_+)$ and $0\leq
t_1<t_2<\cdots<t_{k+1}<T<\infty$, there exists a sequence $h_n$ with
$\lim_{n\to\infty} h_n=0$ and
\[
\sup_{t\in[t_{k+1},T]}\bigl|E\bigl[\Phi\bigl(\check Y^{(n)}_{t_1},\ldots,\check Y^{(n)}_{t_k}\bigr) \bigl(\phi\bigl(\check
Y^{(n)}_{t+h_n}\bigr)-\phi\bigl(\check Y^{(n)}_{t}
\bigr)-h_n\check L\phi\bigl(\check Y^{(n)}_t
\bigr)\bigr)\bigr]\bigr|=o(h_n),
\]
where $\check{\mathcal{L}}$ is given as in (\ref{eq8.2}).

Let for $\phi\in C^\infty_c(\IR_+)$ and $(x,y)\in[1,\infty)\times
\IR_+$
\[
\mathcal{L}^{(n)}_x\phi(y):=\lambda_2^{(n)}
n^2 x y \sum_{k=0}^\infty
\biggl[\phi \biggl(y+\frac{k-1}{n} \biggr)-\phi (y) \biggr]
\mu_2^{(n)}(k)
\]
and recall $\mathcal{L}_{x}$ from (\ref{eq8.6}).
Then
%
\begin{eqnarray}
\label{eqabbranchingcase}\quad
&&
E \bigl[\Phi \bigl(\check Y^{(n)}_{t_ 1},\ldots,\check Y^{(n)}_{t_k} \bigr) \bigl(\phi\bigl(\check
Y^{(n)}_{t+h_n}\bigr)-\phi\bigl(\check Y^{(n)}_{t}
\bigr) \bigr) \bigr]
\nonumber\\
&&\qquad=E \biggl[\Phi \bigl(\check Y^{(n)}_{t_1},\ldots,\check
Y^{(n)}_{t_k} \bigr) \int_t^{t+h_n}
\mathcal{L}^{(n)}_{\check X^{(n)}_s}\phi\bigl(\check Y^{(n)}_t
\bigr)\,ds \biggr]
\\
&&\qquad\quad{} +E \biggl[\Phi \bigl(\check Y^{(n)}_{t_1},\ldots,\check
Y^{(n)}_{t_k} \bigr) \int_t^{t+h_n}
\bigl[ \mathcal{L}^{(n)}_{
\check X^{(n)}_s}\phi \bigl(\check
Y^{(n)}_s \bigr)-\mathcal{L}^{(n)}_{
\check X^{(n)}_s}
\phi \bigl(\check Y^{(n)}_t \bigr) \bigr]\,ds \biggr].
\nonumber
\end{eqnarray}
Using Lemma~\ref{lemsuptsupnexpmoment}, we get, as in (\ref
{eqsecondterm}), that the second term in the last display is $o(h_n)$
uniformly in $t\in[t_{k+1},T]$. Thus
%
\begin{eqnarray}
\label{eqaveragingequation}
&&
E \bigl[\Phi \bigl(\check Y^{(n)}_{t_1},\ldots,\check Y^{(n)}_{t_k} \bigr)
\bigl(\phi\bigl(\check
Y^{(n)}_{t+h_n}\bigr)-\phi\bigl(\check Y^{(n)}_{t}
\bigr)-h_n\check{\mathcal{L}}\phi\bigl(\check Y^{(n)}_t
\bigr) \bigr) \bigr]
\nonumber\\
&&\qquad=E \biggl[\Phi \bigl(\check Y^{(n)}_{t_1},\ldots,\check
Y^{(n)}_{t_k} \bigr) \biggl(\int_{t}^{t+h_n}
\bigl(\mathcal{L}^{(n)}_{ \check
X^{(n)}_s}\phi \bigl(\check
Y^{(n)}_t \bigr)-\mathcal{L}_{\check
X^{(n)}_s}\phi \bigl(
\check Y^{(n)}_t \bigr) \bigr) \,ds \biggr) \biggr]
\nonumber\\[-8pt]\\[-8pt]
&&\qquad\quad{}+ E \biggl[\Phi \bigl(\check Y^{(n)}_{t_1},\ldots,\check
Y^{(n)}_{t_k} \bigr) h_n \biggl(\frac{1}{h_n}
\int_{t}^{t+h_n} \mathcal{L}_{ \check
X^{(n)}_s}\phi
\bigl(\check Y^{(n)}_t\bigr)\,ds-\check{\mathcal{L}}\phi\bigl(
\check Y^{(n)}_t\bigr) \biggr) \biggr]\nonumber\\
&&\qquad\quad{}+o(h_n).
\nonumber
\end{eqnarray}
Calculations similar to those in the proof of Theorem \ref
{thmcrdiffusion} show that the first term on the right-hand side in
the last display is $o(h_n)$ uniformly in $t\in[t_{k+1},T]$ [see proof
of (\ref{eqH4forsingletypecont})], while the second term can be
written as
\begin{eqnarray*}
&&
E \biggl[\Phi \bigl(\check Y^{(n)}_{t_1},\ldots,\check
Y^{(n)}_{t_k} \bigr) h_n \biggl(c_2
\lambda_2\check Y^{(n)}_t\phi'
\bigl(\check Y^{(n)}_t\bigr)+\frac
{1}{2}
\alpha_2\lambda_2 \check Y^{(n)}_t
\phi''\bigl(\check Y^{(n)}_t
\bigr) \biggr)
\\
&&\hspace*{161.7pt}{}\times \biggl(\frac{1}{h_n}\int_{t}^{t+h_n} \check
X^{(n)}_s \,ds - m_X \biggr) \biggr].
\end{eqnarray*}
To show that the latter term is $o(h_n)$ uniformly in $t\in
[t_{k+1},T]$, it suffices
to show the following result.

\begin{theorem}\label{th8.1}
As $n\to\infty$
\[
\sup_{t\in[t_{k+1},T]} E\biggl\llvert E_{\check{\mathcal{F}}^{(n)}_t}
\biggl(\frac{1}{h_n} \int_{t}^{t+h_n}\check X^{(n)}_s \,ds-m_X
\biggr)\biggr\rrvert \to0,
\]
where $\check{\mathcal{F}}^{(n)}_t:=\sigma\{(\check X^{(n)}_s, \check
Y^{(n)}_s)\dvtx s\le t\}$ and $E_{\check{\mathcal{F}}^{(n)}_t}(\cdot) =
E(\cdot|\check{\mathcal{F}}^{(n)}_t)$.
\end{theorem}

In order to prove this theorem, we need the following three results.
Let $S:=D(\IR_+\dvtx [1,\infty)\times\IR_+)$, $\mathcal{P}(S)$ be the
space of probability measures on $S$, and, given a sequence $\{t_n\}
\subset[t_{k+1},T]$, $\mu_n$ be a sequence of $\mathcal{P}(S)$
valued random variables defined as follows. For $A\in\mathcal{B}(S)$,
%
\begin{equation}
\label{occmzr} \mu_n(A)=\frac{1}{a_nh_n}\int_{a_nt_n}^{a_nt_n+a_nh_n}P
\bigl[\bigl(\hat X^{(n)}_{s+\cdot},\hat\eta^{(n)}_{s+\cdot}-
\hat\eta^{(n)}_s\bigr)\in A|\check{\mathcal{F}}^{(n)}_{t_n}\bigr]\,ds.
\end{equation}
Let $S_0:=C(\IR_+\dvtx [1,\infty)\times\IR_+)$.

\begin{lemma}\label{lemrandommeasuretight}
The family of $\mathcal{P}(S)$ valued random variables $\{\mu_n\}_{n\in\IN}$ is tight, and any weak limit point is a $\mathcal
{P}(S_0)$ valued random variable.
\end{lemma}
Let $\pi=(\pi^{(1)},\pi^{(2)})$ with $\pi^{(1)}$ and $\pi^{(2)}$
being the canonical coordinate processes on $S_0$.
%
\begin{lemma}\label{lemrandommeasurelimit}
Let $\mu$ be a weak limit point of $\{\mu_n\}$ given on some
probability space $(\Omega_0,\mathcal{F}_0,P_0)$. Then for $P_0$ a.e.
$\omega\in\Omega_0$, $\mu(\omega)$ satisfies the following:
\begin{enumerate}[(d)]
\item[(a)] $\mu(w)(\pi^{(1)}(t+\cdot)\in F)=\mu(\omega)(\pi^{(1)}\in
F)$, for all $t\ge0$; $F\in\mathcal{B}(C(\IR_+\dvtx [1,\infty)))$;
\item[(b)] $\pi^{(2)}$ is nondecreasing and $\pi^{(2)}_0=0$ a.s.
$\mu(\omega)$;
\item[(c)] $\int_0^\infty1_{(1,\infty)}(\pi^{(1)}_u)\,d\pi^{(2)}_u=0$ a.s. $\mu(\omega)$;
\item[(d)] under $\mu(\omega)$, for all $\phi\in C^\infty_c([1,\infty))$
\[
\phi\bigl(\pi^{(1)}_t\bigr)- \phi\bigl(
\pi^{(1)}_0\bigr)-\int_0^t
\mathcal{L}_1\phi \bigl(\pi^{(1)}_s\bigr)\,ds -
\phi'(1)\pi^{(2)}_t
\]
is a $\{\mathcal{G}_t\}$-martingale, where $\mathcal{L}_1$ is as in
(\ref{eq7.2}) and \mbox{$\mathcal{G}_t:=\sigma\{(\pi^{(1)}_s,\pi^{(2)}_s)\dvtx  s\le t\}$}.
\end{enumerate}
\end{lemma}
We postpone the proofs of Lemmas~\ref{lemrandommeasuretight} and
\ref{lemrandommeasurelimit} until after the proof of Theorem~\ref
{th8.1}. The following is
immediate from the above two lemmas, Proposition \ref
{propstationarydistrofX} and the martingale characterization of the
probability law
of the process in (\ref{eqX}); see Theorem~\ref{thmmartingalecharacterization}.

\begin{corollary}\label{lemrandommeasurelimitcorollary}
Let $(X,\eta)$ be as in Proposition~\ref{propuniqunessofX,Y}
with $X_0\sim\nu_1$ and $\nu_1$ given as in Proposition \ref
{propstationarydistrofX}.
Let
$\mu_0$ be the probability measure on $S_0$ induced by $(X,\eta)$.
Then $\mu_n$ converges weakly to $\mu_0$.
\end{corollary}
\begin{pf*}{Proof of Theorem~\ref{th8.1}}
It suffices to show that for an arbitrary sequence $\{t_n\}\subset
[t_{k+1},T]$ we have, as $n\to\infty$,
\begin{eqnarray*}
&&
E\biggl\llvert E_{\check{\mathcal{F}}^{(n)}_{t_n}} \biggl(\frac{1}{h_n}\int
_{t_n}^{t_n+h_n}\check X^{(n)}_s
\,ds-m_X \biggr)\biggr\rrvert
\\
&&\qquad= E\biggl\llvert E_{\check{\mathcal{F}}^{(n)}_{t_n}} \biggl(\frac
{1}{a_nh_n}\int
_{a_nt_n}^{a_nt_n+a_nh_n}\hat X^{(n)}_s
\,ds-m_X \biggr)\biggr\rrvert \to0.
\end{eqnarray*}
Since
\[
E_{\check{\mathcal{F}}^{(n)}_{t_n}} \biggl(\frac{1}{a_nh_n}\int_{a_nt_n}^{a_nt_n+a_nh_n}
\hat X^{(n)}_s \,ds \biggr)=\int\pi^{(1)}_0\,d
\mu_n,
\]
it suffices to show that
%
\begin{equation}
\label{eq8.14} E\biggl\llvert \int\pi^{(1)}_0\,d
\mu_n - \int\pi^{(1)}_0\,d\mu_0
\biggr\rrvert \to0 \qquad\mbox{as } n\to\infty.
\end{equation}
For any $c>0$, let $\psi_c$ be the following continuous function:
\[
\psi_c(x)= %
\cases{1, &\quad if $x\le\dfrac{c}{2}$,
\vspace*{2pt}\cr
0, &\quad if $x \ge c$,} %
\]
and $\psi_c$ is linearly interpolated on $[\frac{c}{2},c]$.
By Corollary~\ref{lemrandommeasurelimitcorollary} $\mu_n$
converges weakly to~$\mu_0$, and therefore, for every $c>0$,
\[
E\biggl\llvert \int\pi^{(1)}_0\psi_c\bigl(
\pi^{(1)}_0\bigr)\,d\mu_n-\int
\pi^{(1)}_0\psi_c\bigl(\pi^{(1)}_0
\bigr)\,d\mu_0\biggr\rrvert \to0 \qquad\mbox{as }n\to \infty.
\]
Moreover, using the estimate in Lemma~\ref{lemsuptsupnexpmoment},
\begin{eqnarray*}
&&
\sup_{n\in\IN} \biggl(E\biggl\llvert \int\pi^{(1)}_0
\bigl(1-\psi_c\bigl(\pi^{(1)}_0\bigr) \bigr)\,d
\mu_n\biggr\rrvert \biggr)
\\
&&\qquad\le\sup_{n\in\IN} \biggl(\frac{1}{a_nh_n}\int_{a_nt_n}^{a_nt_n+a_nh_n}E
\bigl(\hat X^{(n)}_s 1_{|\hat X^{(n)}_s|\ge
{c}/{2}} \bigr)\,ds \biggr)
\to0 \qquad\mbox{as }c\to\infty
\end{eqnarray*}
and
\[
E\biggl\llvert \int\pi^{(1)}_0\bigl(1-\psi_c
\bigl(\pi^{(1)}_0\bigr)\bigr)\,d\mu_0\biggr
\rrvert \le E ( X_0 1_{|X_0|\ge{c}/{2}} )\to0 \qquad\mbox{as }c\to\infty.
\]
The last three displays imply the convergence in (\ref
{eq8.14}), and thus the result follows.
\end{pf*}
\begin{pf*}{Proof of Lemma~\ref{lemrandommeasuretight}}
To show the tightness of $\{\mu_n\}$, it suffices to show that $\{\nu_n\}$ is tight, where for $A\in\mathcal{B}(S)$
\[
\nu_n(A):=E\mu_n(A)=\frac{1}{a_nh_n}\int
_{a_nt_n}^{a_nt_n+a_nh_n}P\bigl[\bigl(\hat X^{(n)}_{s+\cdot},
\hat\eta^{(n)}_{s+\cdot}-\hat\eta^{(n)}_s
\bigr)\in A\bigr]\,ds.
\]
However, the tightness of $\nu_n$ is immediate in view of the
tightness of
\[
\bigl\{\bigl(\hat X^{(n)}_{s+\cdot},\hat\eta^{(n)}_{s+\cdot
}-\hat\eta^{(n)}_s\bigr)\bigr\}_{n\in\IN,s\in\IR_+},
\]
which was proved in Proposition~\ref{lemtighnessofcatalyst-reactant}.

Let $\mu$ be a weak limit point of $\mu_n$ and $J\dvtx S\to\IR_+$ be
defined by
\[
J(\pi):=\int_0^\infty e^{-u}
\bigl[J(\pi,u)\wedge1\bigr]\,du,
\]
where
\[
J(\pi,u):=\sup_{0\le t\le u} \bigl(\bigl|\Delta\bigl(\pi^{(1)}_t
\bigr)\bigr|+\bigl|\Delta \pi^{(2)}_t\bigr| \bigr).
\]
Then $J$ is continuous and bounded on $S$, and in order to show that
$\mu$ is supported on $S_0$, it suffices to show that
$
\mu(J(\pi)=0)=1
$; see~\cite{Ethier1986}, page 147.
In turn, for the latter equality to hold, it suffices to show that for
all $\varepsilon>0$,
$
E\mu_n(J(\pi)>\varepsilon)\to0,
$
as $n\to\infty$.
Now
\begin{eqnarray*}
&&
E\mu_n\bigl(J(\pi)>\varepsilon\bigr)
\\
&&\qquad=\frac{1}{a_nh_n} 
\int_{a_nt_n}^{a_nt_n+a_nh_n} P \biggl(
\int_0^\infty e^{-u} \Bigl(
\sup_{ t\le u}\bigl[ \bigl|\Delta\hat X^{(n)}_{s+t}\bigr|\\
&&\qquad\quad\hspace*{151.5pt}{}+\bigl|\Delta
\bigl(\hat \eta^{(n)}_{s+t}-\hat\eta^{(n)}_s
\bigr)\bigr|\bigr]\wedge1 \Bigr)\,du>\varepsilon \biggr) \,ds.
\end{eqnarray*}
Finally, noting that $\hat\eta^{(n)}_{s+\cdot}-\hat\eta^{(n)}_s$
is continuous and using Lemma~\ref{lemPsupDeltahatX^nt}, we now
have that the right-hand side of the latter equation converges to 0 as
$n\to\infty$. The result follows.
\end{pf*}
\begin{pf*}{Proof of Lemma~\ref{lemrandommeasurelimit}}
For a measure $\nu\in\mathcal{P}(S)$, let $\tilde E^{\nu}$ denote
the expectation operator.
For (a), we show that
%
\begin{eqnarray}
\label{eq8.15}
&&
\tilde E^{\mu(\omega)}\bigl( f\bigl(\pi^{(1)}_{t+\cdot}
\bigr)\bigr)-\tilde
E^{\mu(\omega)}\bigl(f\bigl(\pi^{(1)}\bigr)\bigr)=0\nonumber\\[-8pt]\\[-8pt]
&&\eqntext{\mbox{a.s. for all bounded continuous } f\mbox{ on }S.}
\end{eqnarray}
Note that
\begin{eqnarray*}
&&\bigl\llvert \tilde E^{\mu_n}\bigl(f\bigl(\pi^{(1)}_{t+\cdot}
\bigr)\bigr)-\tilde E^{\mu
_n}\bigl(f\bigl(\pi^{(1)}\bigr)\bigr)
\bigr\rrvert
\\
&&\qquad =\biggl\llvert \frac{1}{a_nh_n}\int_{a_nt_n}^{a_nt_n+a_nh_n}
\bigl(E_{\check{\mathcal{F}}^{(n)}_{t_n}} f\bigl(\hat X^{(n)}_{s+t+\cdot}\bigr) -
E_{\check{\mathcal{F}}^{(n)}_{t_n}} f\bigl(\hat X^{(n)}_{s+\cdot}\bigr) \bigr) \,ds
\biggr\rrvert
\\
&&\qquad \le\frac{2t}{a_nh_n}\|f\|_{\sup}\to0 \qquad\mbox{as }n\to \infty.
\end{eqnarray*}
This proves (\ref{eq8.15}) since we can choose
$h_n$ such that $a_n h_n\to\infty$, and thus (a) follows.

Property (b) is immediate from the fact that $\hat\eta^{(n)}_{s+\cdot
}-\hat\eta^{(n)}_s$ is nondecreasing and continuous with initial
value 0 for each $n$.\vadjust{\goodbreak}

To prove (c), it suffices to show that for a.e. $\omega$ and for
every $T,\delta>0$
\[
\tilde E^{\mu(\omega)} \biggl(\int_0^T
f_\delta\bigl(\pi^{(1)}_s\bigr)\,d
\pi^{(2)}_s\wedge1 \biggr)=0,
\]
where $f_\delta$ is defined in (\ref{eqfdelta}). In turn, for the
above equality to hold, it suffices to show that
\[
E \biggl(\tilde E^{\mu} \biggl(\int_0^T
f_\delta\bigl(\pi^{(1)}_s\bigr)\,d
\pi^{(2)}_s\wedge1 \biggr) \biggr)=0.
\]
The latter equality is immediate on noting that for every $T,\delta>0$
\begin{eqnarray*}
&&
E \biggl(\tilde E^{\mu_n} \biggl(\int_0^T
f_\delta\bigl(\pi^{(1)}_u\bigr)\,d
\pi^{(2)}_u\wedge1 \biggr) \biggr)
\\[-2pt]
&&\qquad =\frac{1}{a_nh_n}\int_{a_nt_n}^{a_nt_n+a_nh_n}E \biggl(\int
_0^T f_\delta\bigl(\hat
X^{(n)}_{s+u}\bigr)\,d\bigl(\hat\eta^{(n)}_{s+\cdot}-
\hat \eta^{(n)}_s\bigr) (u)\wedge1 \biggr) \,ds=0
\end{eqnarray*}
and thus
\[
E \biggl[\tilde E^{\mu} \biggl(\int_0^T
f_\delta\bigl(\pi^{(1)}_s\bigr)\,d
\pi^{(2)}_s\wedge1 \biggr) \biggr]=\lim_{n\to\infty} E
\biggl[\tilde E^{\mu_n} \biggl(\int_0^T
f_\delta\bigl(\pi^{(1)}_s\bigr)\,d
\pi^{(2)}_s\wedge1 \biggr) \biggr]=0.
\]
Finally, consider (d).
It suffices to show that for every $0\le r\le t<\infty$,
%
\begin{eqnarray*}
&&
E \biggl| \tilde E^{\mu} \biggl( \psi\bigl(\pi^{(1)},
\pi^{(2)}\bigr) \biggl(\phi \bigl(\pi^{(1)}_t\bigr)-
\phi\bigl(\pi^{(1)}_r\bigr) -\int_r^t
\mathcal{L}_1\phi\bigl(\pi^{(1)}_u\bigr)\,du
\\[-2pt]
&&\hspace*{151.5pt}\qquad{}- \phi'(1)\bigl[\pi_t^{(2)}-
\pi_r^{(2)}\bigr] \biggr) \biggr) \biggr|=0,
\end{eqnarray*}
where $\psi\dvtx S\to\IR$ is an arbitrary bounded, continuous, $\mathcal
{G}_s$ measurable map.

Now fix such $r,t$ and $\psi$. Assume without loss of generality that
$\mu_n$ converges to $\mu$. Combining this weak convergence with
Lemma~\ref{lemEVofsquaredsupremaofshiftedprocesses}, we see
that the left-hand side of the last display is the limit of
\begin{eqnarray*}
&& E \biggl|\tilde E^{\mu_n} \biggl( \psi\bigl(\pi^{(1)},
\pi^{(2)}\bigr) \biggl(\phi\bigl(\pi^{(1)}_t\bigr)-
\phi\bigl(\pi^{(1)}_r\bigr) -\int_r^t
\mathcal{L}_1\phi \bigl(\pi^{(1)}_u\bigr)\,du\\[-2pt]
&&\qquad\quad\hspace*{144.5pt}{}-
\phi'(1)\bigl[\pi_t^{(2)}-\pi_r^{(2)}
\bigr]\biggr) \biggr) \biggr|
\\[-2pt]
&&\qquad=\frac{1}{a_nh_n}\int_{a_nt_n}^{a_nt_n+a_nh_n} E \biggl|
E_{\check{\mathcal{F}}_{t_n}} \biggl(\psi\bigl(\hat X^{(n)}_{s+\cdot},\hat
\eta^{(n)}_{s+\cdot}-\hat\eta^{(n)}_s\bigr)
\\[-2pt]
&&\hspace*{146.2pt}{}\times \biggl[\phi\bigl(\hat X^{(n)}_{s+t}\bigr)-\phi\bigl(\hat
X^{(n)}_{s+r}\bigr) \\[-2pt]
&&\qquad\quad\hspace*{130pt}{}-\int_r^t
\mathcal{L}_1\phi\bigl(\hat X^{(n)}_{s+u}
\bigr)\,du\\[-2pt]
&&\qquad\quad\hspace*{130pt}{}-\phi'(1)\bigl[\hat\eta^{(n)}_{s+t}-\hat
\eta^{(n)}_{s+r}\bigr] \biggr] \biggr) \biggr|\,ds.
\end{eqnarray*}
To complete the proof, it suffices to show that the limit of the
expression in the last display is 0. Note that for
$\phi\in C^\infty_c([1,\infty))$,
\[
\phi\bigl(\hat X_t^{(n)}\bigr)-\phi\bigl(\hat
X_0^{(n)}\bigr)-\int_0^t
\mathcal {L}^{(n)}_1\phi\bigl(\hat X_s^{(n)}
\bigr)\,ds-\mathcal{D}^{(n)}_1\phi(1)\hat
\eta^{(n)}_t
\]
is a martingale, where $\mathcal{D}^{(n)}_1\phi(1):=n[\phi(1)-\phi
(1-\frac{1}{n})]$ and
\[
\mathcal{L}^{(n)}_1\phi(x):=\lambda_1^{(n)}
n^2 x\sum_{k=0}^\infty \biggl[
\phi \biggl(x+\frac{k-1}{n} \biggr)-\phi(x) \biggr]\mu_1^{(n)}(k).
\]
Thus, it suffices to prove that
\[
\lim_{n\to\infty}\frac{1}{a_nh_n}\int_{a_nt_n}^{a_nt_n+a_nh_n}E
\biggl\llvert \int_r^t \bigl(
\mathcal{L}^{(n)}_1\phi\bigl(\hat X_{s+u}^{(n)}
\bigr)-\mathcal {L}_1\phi\bigl(\hat X_{s+u}^{(n)}
\bigr) \bigr)\,du\biggr\rrvert \,ds=0
\]
and
\[
\lim_{n\to\infty}\frac{1}{a_nh_n}\int_{a_nt_n}^{a_nt_n+a_nh_n}
\bigl\llvert \mathcal{D}^{(n)}_1 \phi(1)-
\phi'(1) \bigr\rrvert E \bigl(\bigl\llvert \hat\eta^{(n)}_{s+t}-
\hat\eta^{(n)}_{s+r}\bigr\rrvert \bigr) \,ds=0.
\]
The proofs for the last two equalities are completed as those for
(\ref{eqH4forsingletypecont})
and
(\ref{eqEintD^n-phi}) upon using the uniform estimates in
Corollary~\ref{corEVofsquaredsupremaofshiftedhatX^n}
and Lemma~\ref{lemEVofsquaredsupremaofshiftedprocesses}.
\end{pf*}

\begin{appendix}\label{app}
\section*{Appendix}
\begin{pf*}{Proof of Proposition~\ref{propuniqunessofX,Y}}
We will consider here only the case where $(X_0, Y_0) \equiv(x,y)$ for
some $(x,y) \in[1, \infty) \times[0, \infty)$. The general case can
be treated similarly.
The unique solvability of (\ref{eqX}) is an immediate consequence of
the Lipschitz property of the Skorohod map, Lipschitz coefficients
(note that $f(x)=\sqrt{x}$ is a Lipschitz function on $[1,\infty)$)
and a standard Picard iteration scheme; see, for example, Proposition 1
in~\cite{AnOr}.

We next argue the unique solvability of (\ref{eqY}). For $n\in\IN$,
let $\sigma^{(n)}:=\inf\{t>0|X_t\ge n\}$, $\bar X^{(n)}_t:=
X_{t\wedge\sigma^{(n)}}$ and $f^{(n)}(y):=y\vee\frac{1}{n}$.
Consider the equation
%
\setcounter{equation}{0}
\begin{eqnarray}
\label{eqcheckYinpropuniqunessofX,Y} \bar Y^{(n)}_t&=&
Y_0+c_2\lambda_2\int_0^t
\bar X^{(n)}_sf^{(n)}\bigl( \bar
Y^{(n)}_s\bigr) \,ds\nonumber\\[-8pt]\\[-8pt]
&&{}+\sqrt{\alpha_2
\lambda_2}\int_0^t \sqrt{
\bar X^{(n)}_s f^{(n)}\bigl( \bar
Y^{(n)}_s\bigr)} \,dB_s^Y.\nonumber
\end{eqnarray}
From the Lipschitz property of $f^{(n)}$ and $\sqrt{f^{(n)}}$ it
follows that, for each $n$, the above equation has a unique pathwise
solution. Let $\tau^{(n)}:=\inf\{t>0|\bar Y^{(n)}_t=\frac{1}{n}\}$
and $\theta^{(n)}:=\tau^{(n)}\wedge\sigma^{(n)}$. Note that $\bar
Y^{(n)}$ solves (\ref{eqY}) on $[0,\theta^{(n)}]$. Also, by unique
solvability of (\ref{eqcheckYinpropuniqunessofX,Y}), we have
for all $n\in\IN$,
$
\bar Y^{(n+1)}(\cdot\wedge\theta^{(n)})=\bar Y^{(n)}(\cdot\wedge
\theta^{(n)}).
$
Finally, letting $\theta^{(\infty)}:=\lim_{n\to\infty}\theta^{(n)}$,
the unique
solution of (\ref{eqY}) is given by the following:
\[
Y_t(\omega) = %
\cases{\bar Y^{(n)}_t(
\omega), &\quad if $0\leq t\leq\theta^{(n)}(\omega )$ for some $n\in\IN$,
\vspace*{2pt}\cr
0, &\quad if $t\geq\theta^\infty(\omega)$.} %
\]
\upqed\end{pf*}
\end{appendix}

\section*{Acknowledgments}

We thank the Editor, the Associate Editor and the referees for their
comments that led to several improvements in the manuscript.



\printaddresses

\end{document}